\documentclass[letter, 12pt, oneside,reqno]{amsart}

\usepackage{latexsym,amsxtra,amscd,ifthen}
\usepackage{amsfonts}
\usepackage{verbatim}
\usepackage{dsfont}
\usepackage{mathtools}
\usepackage{enumerate}
\usepackage[capitalise]{cleveref}
\usepackage{endnotes}
\usepackage{graphicx}
\usepackage[latin1]{inputenc}
\usepackage{array}
\usepackage{amsmath}
\usepackage{fancyhdr}
\usepackage{graphics}
\usepackage{amssymb}
\usepackage{amsthm}
\usepackage{enumitem}
\usepackage{mathrsfs}
\usepackage{yfonts} 
\usepackage{cancel}
\usepackage[square,numbers]{natbib}
\usepackage{stmaryrd}
\usepackage{amsfonts} 
\usepackage[pdftex]{color}
\usepackage{ulem}
\usepackage{tikz-cd}		
\usepackage[arrow, matrix, curve]{xy}

\newlength{\sublength}
\newlength{\suplength}
\newcommand\lsubsup[2]{%
  \settowidth{\sublength}{#1}%
  \settowidth{\suplength}{#2}%
  _{\text{%
      \ifdim \sublength<\suplength
        \hphantom{#2}\llap{#1}%
      \else
        #1%
      \fi
    }%
  }%
  ^{\text{%
      \ifdim \suplength<\sublength
        \hphantom{#1}\llap{#2}%
      \else
        #2%
      \fi
    }%
  }%
}

\makeindex

\AtBeginDocument{
  \label{CorrectFirstPageLabel}
  
}

\begin{document}

\def\mathllap{\mathpalette\mathllapinternal}
\def\mathllapinternal#1#2{%
\llap{$\mathsurround=0pt#1{#2}$}
}
\def\clap#1{\hbox to 0pt{\hss#1\hss}}
\def\mathclap{\mathpalette\mathclapinternal}
\def\mathclapinternal#1#2{%
\clap{$\mathsurround=0pt#1{#2}$}%
}
\def\mathrlap{\mathpalette\mathrlapinternal}
\def\mathrlapinternal#1#2{%
\rlap{$\mathsurround=0pt#1{#2}$}
}

\def\itemMath#1{\raisebox{-\abovedisplayshortskip}
{\parbox{1.0\linewidth}{\begin{equation}\begin{split}#1\end{split}\end{equation}}}}

\def\itemMathtwo#1#2{\raisebox{-#1pt}
{\parbox{1.0\linewidth}{\begin{equation}\begin{split}#2\end{split}\end{equation}}}}

\newenvironment{other}[1]{\noindent\textbf{\stepcounter{thm}\thethm\ #1}\hspace*{1,5em}}{}

\let\amsamp=&

\newcommand{\CC}{\mathbb{C}}
\newcommand{\NN}{\mathbb{N}}
\newcommand{\QQ}{\mathbb{Q}}
\newcommand{\ZZ}{\mathbb{Z}}
\newcommand{\RR}{\mathbb{R}}

\newcommand{\R}{\Rightarrow}
\newcommand{\id}{\mathrm{id}}
\newcommand{\Sym}{\operatorname{Sym}}
\newcommand{\Aut}{\operatorname{Aut}}
\newcommand{\End}{\operatorname{End}}
\newcommand{\Hom}{\operatorname{Hom}}

\newcommand{\HS}{\operatorname{H}}
\newcommand{\Tr}{\operatorname{Tr}}
\newcommand{\TrV}{\underline{\operatorname{Tr}}}
\newcommand{\TrM}{\underline{\operatorname{Tr_M}}}
\newcommand{\HV}{\underline{\operatorname{H}}}
\newcommand{\HM}{\underline{\operatorname{H_M}}}
\newcommand{\tr}{\operatorname{tr}}
\newcommand{\gr}{\operatorname{gr}}
\newcommand{\op}{\operatorname{op}}
\newcommand{\hdet}{\operatorname{hdet}}
\newcommand{\Hdet}{\operatorname{Hdet}}
\newcommand{\grA}{\operatorname{gr}-\operatorname{A}}
\newcommand{\Ext}{\operatorname{Ext}}
\newcommand{\GK}{\operatorname{GKdim}}
\newcommand{\GKt}{\operatorname{GK}}
\newcommand{\gldim}{\operatorname{gl.dim}}
\newcommand{\pdim}{\operatorname{pdim}}
\newcommand{\injdim}{\operatorname{injdim}}
\newcommand{\proj}{\operatorname{proj}\,}
\newcommand{\cd}{\operatorname{cd}}
\newcommand{\qgr}{\operatorname{qgr}}
\newcommand{\tors}{\operatorname{tors}}
\newcommand{\sgn}{\operatorname{sgn}}
\newcommand{\type}{\operatorname{type}}
\newcommand{\Tor}{\operatorname{Tor}}
\newcommand{\adj}{\operatorname{adj}}
\newcommand{\col}{\operatorname{col}}

\newtheoremstyle{paper}{15pt}{15pt}{}{}{\bf}{}{1.5em}{}
\newtheoremstyle{proof}{15pt}{15pt}{}{}{\bfseries\itshape}{}{\newline}{}

\theoremstyle{paper}
\newtheorem{thm}{Theorem}[section]
\newtheorem*{thm2}{Theorem}
\newtheorem{de}[thm]{Definition}
\newtheorem{cor}[thm]{Corollary}
\newtheorem{rem}[thm]{Remark}
\newtheorem{lemma}[thm]{Lemma}
\newtheorem{ex}[thm]{Example}
\newtheorem{prop}[thm]{Proposition}

\theoremstyle{proof}
\newtheorem*{pro}{Proof}
\newtheorem*{proof2}{Idea of Proof}

\newcommand*\alt{}
\let\alt\thethm
\renewcommand\thethm{(\alt)}

\title{Properties of the fixed ring of a preprojective algebra}
\author{Stephan Weispfenning\protect\endnotemark[1]}



\address{1. To whom correspondence should be addressed\hfill\break E-mail: sweispfe@ucsd.edu\hfill\break The author declares no conflict of interest.}
	
\begin{abstract}
For a finite group acting on a polynomial ring, the Chevalley-Shephard-Todd Theorem proves that the fixed subring is isomorphic to a polynomial ring if and only if the group is generated by pseudo-reflections. In recent years, progress was made in work of Kirkman, Kuzmanovich, Zhang, and others to extend this result to Artin-Schelter regular algebras by expanding pseudo-reflections to quasi-reflections. Naturally, the question arises if the theory generalizes further to non-connected noncommutative algebras. Our objects of study will be preprojective algebras which are certain factor algebras of path algebras corresponding to extended Dynkin diagrams of type $A$, $D$ or $E$. This work answers the question what conditions need to be satisfied by the fixed ring in order to make a rich theory possible. On our way, we will point out additional difficulties in establishing quasi-reflections using the trace and reveal situations which do not occur for Artin-Schelter regular algebras.
\end{abstract}

\maketitle

\pagestyle{fancy}
\fancyhf{}
\pagenumbering{arabic}\cfoot{-\ \thepage\ -}

\setlength{\headheight}{17pt}

\fancyhead[L]{\ifthenelse{\isodd{\value{page}}}{\center\small \small\rightmark}{\center STEPHAN WEISPFENNING}}

\section{Introduction}
In 1954 G.~C.~Shephard and J.~A.~Todd proved one of the most powerful theorems in commutative invariant theory. Shortly after, C.~Chevalley provided a uniform proof in 1955. Now known as the Chevalley-Shephard-Todd Theorem, it says that for a finite group $G$ acting faithfully on $k[x_1,\ldots,x_n]$ its fixed ring $k[x_1,\ldots,x_n]^G$ is isomorphic to $k[x_1,\ldots,x_n]$ if and only if $G$ is generated by pseudo-reflections. An element of $G$ is a pseudo-reflection if its restriction to the degree $1$ piece of $k[x_1,\ldots,x_n]$ fixes a hyperplane. In 1987 M.~Artin and W.~Schelter introduced a new class of noncommutative graded algebras which can be understood as noncommutative generalizations of polynomial rings: an Artin-Schelter regular algebra is defined to be a finitely generated connected $\NN$-graded ring $S = k \oplus \bigoplus_{j > 0}{S_j}$ satisfying the following three conditions: \vspace{-6pt}
\begin{itemize}[topsep=0pt, partopsep=0pt,leftmargin=0.75in]
	\item[(AS 1)] $S$ has finite (graded) global dimension $d$. \par
	\item[(AS 2)] $S$ has finite Gelfand-Kirillov dimension. \par
	\item[(AS 3)] $\Ext_S^m(k,S) \cong \Ext_{S^{\op}}^m(k,S) \cong \delta_{md}k[\ell]$ for some $\ell$. \par
\end{itemize}

Chan, J{\o}rgensen, Kirkman, Kuzmanovich, Walton, Zhang and many others generalized the theory to finite group actions (\cite{MR2434290, MR2601538, MR1758250}) and finite Hopf algebra actions (\cite{MR2568355, MR3552496}) on Artin-Schelter regular algebras. In some of their works a more general version of the Chevalley-Shephard-Todd Theorem was found. 

After expanding the polynomial ring to an Artin-Schelter regular algebra $S$, the first question was to find the suitable generalization of pseudo-reflections. The definition to have a pseudo-reflection fix a hyperplane of the degree $1$ piece fails terribly in this new setting (see \cite[Example 1.4]{MR3546709}). Instead, Jing and Zhang showed in \cite[Theorem 2.3]{MR1438180} that the formal power series, called the trace of $g$,
\[
	\Tr_S(g,t) = \sum_{j \in \NN}{\tr(g|_{S_j})t^j}
\]
has the special form $e_g(t)^{-1}$ for some polynomial $e_g(t) \in k[t]$. Moreover, for a nontrivial automorphism $g$ acting on an Artin-Schelter regular domain $S$, the multiplicity of $1$ as a root of $e_g(t)$ is always smaller than $\GK(S)$ (\cite[Proposition 1.8]{MR2434290}). The term quasi-reflection was introduced in \cite[Definition 2.2]{MR2434290} to describe graded automorphisms $g$ for which $1$ is a root of $e_g(t)$ of multiplicity $\GK(S) - 1$. This definition coincides with the definition of pseudo-reflections for polynomial rings. For quantum polynomial rings, Artin-Schelter regular algebras with Hilbert series equal to the Hilbert series of a polynomial ring with equal $\GKt$-dimension, Kirkman, Kuzmanovich and Zhang classify all quasi-reflections in \cite[Proposition 4.3]{MR2434290}. It turns out that in addition to the standard pseudo-reflections a new class called mystic reflections arises (\cite[Definition 3.9]{MR2434290}). Using the classification of quasi-reflections, the foremost generalization of the Chevalley-Shephard-Todd Theorem was proved by Kirkman, Kuzmanovich and Zhang in \cite{MR2601538}: 

\begin{thm2}[{\cite[Theorem 5.5 and Theorem 6.3]{MR2601538}}]
Let $S$ be a noetherian and Koszul AS-regular domain and $G$ a finite subgroup of $\Aut_{\gr}(S)$. Suppose either $G$ is abelian or $S$ is a skew polynomial ring. Then the fixed subring $S^G$ is Artin-Schelter regular if and only if $G$ is generated by quasi-reflections.
\end{thm2}
Additionally, the proof of this beautiful result uses the observation that if the fixed ring $S^G$ has finite global dimension it is already Artin-Schelter regular (\cite[Lemma 1.10(c)]{MR2434290}). At the same time, this is equivalent to $S$ being a free module over the fixed ring (\cite[Lemma 1.11]{MR2434290}). 

Starting with this paper, the main goal is to find a version of the Chevalley-Shephard-Todd Theorem for a class of non-connected noncommutative algebras. We begin with a quiver $Q$ which for us is a directed graph with finitely many vertices $Q_0 = \{e_i \mid 1 \leq i \leq n\}$ and finitely many arrows $Q_1$. The double $\bar{Q}$ of $Q$ is the quiver obtained by keeping the vertex set and adding a new arrow $\alpha^\ast$ from $e_j$ to $e_i$ for each arrow $\alpha$ from $e_i$ to $e_j$. The preprojective algebra $\Pi_Q(k)$ is then defined as 
\[
	A = \Pi_Q(k) = k\bar{Q}/\left(\sum_{\alpha \in Q_1}{\alpha \alpha^\ast} - \sum_{\alpha \in Q_1}{\alpha^\ast \alpha}\right)
\]
where $k\bar{Q}$ denotes the path algebra of $\bar{Q}$.

Primarily, we focus on preprojective algebras $A$ for the extended Dynkin diagram $\widetilde{A_{n-1}}$. However, the main theory developed in Sections \ref{sectiongeneral} and \ref{section6} applies to noetherian $\NN$-graded twisted Calabi-Yau algebras with degree zero piece $B_0 \cong k^n$. We define twisted Calabi-Yau algebras in Section \ref{section2}, Definition \ref{CYdef} which form a non-connected generalization of Artin-Schelter regular algebras. In fact, Reyes, Rogalski and Zhang showed in \cite[Lemma 1.2]{MR3250287} that every connected twisted Calabi-Yau algebra is indeed Artin-Schelter regular. Preprojective algebras form a nice class of non-connected noncommutative algebras satisfying the mentioned properties. More precisely, they have a natural grading, global dimension and Gelfand-Kirillov dimension equal to $2$ and nice homological properties including the $\chi$-condition. In particular, this is a class of Calabi-Yau algebras (Proposition \ref{prop210}).

The first steps towards finding a Chevalley-Shephard-Todd Theorem consist of finding the right definition of a quasi-reflection and the right condition to require the fixed ring to have. As mentioned earlier, in \cite{MR2601538}'s version a key argument is that finite global dimension already implies the important homological properties to make the fixed ring Artin-Schelter regular (\cite[Lemma 1.10]{MR2434290}). Unfortunately, the situation with non-connected algebras is significantly more difficult. To obtain fixed rings of preprojective algebras which are twisted Calabi-Yau it is not enough to require only finite global dimension. We moreover need to ask for global dimension equal to $\gldim(A)$ (Example \ref{exp43}). Additionally, if the fixed ring is twisted Calabi-Yau the strongest statement possible is that $A$ is projective as an $A^G$-module rather than free (Example \ref{exp42}). The following is the main result of this paper stated in its most general form which applies to preprojective algebras. Furthermore, it answers one of the two main questions before attempting to find a Chevalley-Shephard-Todd Theorem for this class of non-connected noncommutative algebras.

\setcounter{section}{6}
\setcounter{thm}{4}
\begin{thm}
Assume $B$ is a noetherian $\NN$-graded twisted Calabi-Yau algebra of global dimension $d \geq 1$ with degree zero piece $B_0 \cong k^n$. Denote the pairwise orthogonal idempotents by $\{e_1,\ldots,e_n\}$. Let $G$ be a finite subgroup of $\Aut_{\gr}(B)$ such that $g(e_i) = e_i$ for all $g \in G$ and all $i = 1,\ldots,n$. Suppose $\pdim(S_i) = \pdim(S_j)$ and $\pdim(S_i^\ast) = \pdim(S_j^\ast)$ for the simple $B^G$-modules $S_i = e_iB^G/(e_iB^G)_{\geq 1}$ and $S_i^\ast = B^Ge_i/(B^Ge_i)_{\geq 1}$ and all $i,j = 1,\ldots,n$. Then the following are equivalent:
\begin{itemize}[topsep=0pt]
	\item[(1)]	$\gldim(B^G) = \gldim(B)$.
	\item[(2)]	$B^G$ is twisted Calabi-Yau.
	\item[(3)]	$B$ is a projective $B^G$-module.
\end{itemize}
\end{thm}

As a consequence, the result can be stated more elegantly for preprojective algebras $\Pi_Q(k)$ for $Q = \widetilde{A_{n-1}}$ which is done in Section \ref{section7}:

\setcounter{section}{7}
\setcounter{thm}{6}
\begin{thm}
For type $Q = \widetilde{A_{n-1}}$, $A = \Pi_Q(k)$ as defined above and $G$ a finite subgroup of the graded automorphism group of $A$ such that $g(e_i) = e_i$ for all $g \in G$ and all $i = 1,\ldots,n$, the following conditions are equivalent:
\begin{itemize}[topsep=0pt, partopsep=0pt]
	\item[(1)]	$\gldim(A^G) = \gldim(A)$ (i.e. $\gldim(A^G) = 2$).
	\item[(2)]	$A^G$ is twisted Calabi-Yau.
	\item[(3)]	$A$ is a projective $A^G$-module.
\end{itemize}
\end{thm}

One crucial condition in both theorems is to require the automorphisms to fix the idempotents $e_1,\ldots,e_n$. Mainly, there are two reasons for this restriction. Part of the proof that the fixed ring is twisted Calabi-Yau consists of proving the generalized Gorenstein condition. This requires knowing the simple modules of the fixed ring which are in one-to-one correspondence with its pairwise orthogonal idempotents. Only if all automorphisms fix the idempotents it is possible to know the set of simple modules of the fixed ring. Secondly, we would like to use the trace function as defined above which forms one of the crucial tools of invariant theory. It is easy to see that already for $Q = \widetilde{A_{2}}$ there exist graded automorphisms of $\Pi_Q(k)$ which permute the idempotents and whose trace functions vanish entirely. This makes it impossible to gain any information and therefore these automorphisms have to be excluded. Nevertheless, Proposition \ref{prop36} proves that the graded automorphisms which fix the idempotents form a subgroup of finite index in the full group of graded automorphisms.

The other big question for extending the theory to preprojective algebras is to find the correct definition for a quasi-reflection. The trace function has been used to define quasi-reflections before (\cite[Definition 2.2]{MR2601538}) and by forcing our automorphisms to fix the idempotents we expect to gain useful information from it. Moreover, for non-connected algebras and automorphisms which fix the idempotents $e_i$, there is a different way to encode more information. More precisely, the trace function for each indecomposable projective piece $e_iA$ for $i = 1,\ldots,n$ put into a vector becomes valuable. In Section \ref{sectionAn}, we find two concrete formulas to calculate this vector trace of a graded automorphism which fixes the idempotents in case $\widetilde{A_{n-1}}$ (Propositions \ref{prop34} and \ref{prop35}). Since $A = \bigoplus_{i=1}^n {e_iA}$, simply summing up the vector trace yields the traditional trace used before by Jing, Kirkman, Kuzmanovich, Zhang and others. Therefore, these formulas help working out three longer and interesting examples in Section \ref{expchapter} all of which take $Q = \widetilde{A_2}$.

The first example shows that it is possible to have a nice fixed ring such that $A$ is free over $A^G$. However, the second example presents a new situation. Even though the fixed ring has global dimension $2$, $A$ is projective and not free over its fixed ring. Nevertheless, the fixed ring $A^G$ is twisted Calabi-Yau (Example \ref{exp42}). Lastly, a more concerning and fascinating example is worked out in \ref{exp43}. The fixed ring has global dimension $3$ unlike the preprojective algebra which has global dimension $2$. This situation cannot happen for Artin-Schelter regular algebras. The crucial problem is that rather than having only one simple module, the preprojective algebra of Example \ref{exp43} and its fixed ring have three simple modules. Therefore, it happens in this example that two of the simple modules have projective dimension $2$ while the third simple module has projective dimension $3$. This makes a homological property such as the generalized Gorenstein condition impossible. Regardless, the global dimension as the maximum of the projective dimensions of the simple modules stays finite. In light of Theorem \ref{thm77}, it also follows that $A$ is not projective over $A^G$ and so this example should not be included in a Chevalley-Shephard-Todd Theorem. 

When looking at the three examples one notices that the usual trace function is no longer of the form $p(t)^{-1}$ for some polynomial $p(t)$. Instead, it is possible to get rational functions which make a straightforward generalization of quasi-reflection by using the multiplicity of $1$ as a root of the denominator not well-defined (Proposition \ref{prop34}). Simply using the multiplicity of the pole $1$ of the trace function does not solve the problem either. Example \ref{exp43} provides a graded automorphism $g$ whose multiplicity of $1$ as a pole of the trace equals $1 = \GK(A) - 1$. But as described, the corresponding fixed ring $A^{\langle g \rangle}$ has global dimension $3$ and does not satisfy the desired homological properties. Therefore, $g$ cannot be considered a quasi-reflection. Additionally, it is possible to have a nonidentity automorphism whose trace has $1$ as a pole of multiplicity $\GK(A) = 2$. How to define quasi-reflections remains an unsolved problem.

In future work, the author hopes to make progress towards a Chevalley-Shephard-Todd Theorem for preprojective algebras. This paper should be understood as a first evaluation of the situation, solving some technical issues, revealing technical problems and situations we do and do not want to consider as well as finding the right analogue on the algebra side of the Chevalley-Shephard-Todd Theorem.

\section*{Acknowledgment}
Most importantly, the author genuinely wishes to thank Daniel Rogalski for proposing this interesting research topic as well as reading and improving drafts of this paper. The results would not exist without his valuable insight and helpful discussions. The author wishes to thank James J. Zhang for helpful conversations. The author was partially supported by the NSF grants DMS-1201572 and DMS-1601920 as well as the NSA grant H98230-15-1-0317. The author also wishes to thank Efim Zelmanov for his generous financial support that allowed him to fully focus on this project. Additionally, the author is grateful to the referee for numerous valuable and valid comments which improved this paper.

\setcounter{section}{1}
\setcounter{thm}{0}
\section{Preliminaries}\label{section2}
Throughout let $k$ be an algebraically closed field of characteristic $0$. Unless said otherwise, every algebra is $\NN$-graded and assumed to be locally finite, i.e. each graded piece is finite dimensional over $k$. All quivers have finitely many vertices and finitely many arrows. The symbol $A$ is reserved to be the preprojective algebra $\Pi_Q(k)$ of an extended Dynkin quiver $Q$ (see Euclidean graphs in \cite[Chapter VII.2]{MR2335985}) as defined in Definition \ref{def23}. 

The first definitions are necessary to define preprojective algebras. Afterwards, we establish a list of conventions in Definition \ref{def24}. It is worthwhile to point out the $\chi$-condition (Definition \ref{defchi}) and the definition of a twisted Calabi-Yau algebra (Definition \ref{CYdef}) separately. Before Proposition \ref{prop210} collects essential properties of preprojective algebras we state an important result giving equivalent conditions for an algebra to be twisted Calabi-Yau. Important tools for any sort of invariant theory are the construction of a Hilbert series and a trace function which we will extend to a vector and matrix notation starting with Definition \ref{def211}. Moreover, Molien's Theorem (see Lemma \ref{lemma214}) connects the Hilbert series of a fixed ring to trace functions. This background section is finished off by quoting two very general results. The first guarantees that the ring is finitely generated over its fixed ring, and the second makes sure that the fixed ring satisfies the $\chi$-condition whenever the ring has this property.

\begin{de}
Following \cite[Chapter II.1]{MR2197389}, a \textit{quiver} $Q = (Q_0,Q_1,s,t)$ is a quadruple consisting of vertices $Q_0 = \{e_i \mid i \in I\}$, arrows $Q_1$ and two maps $s,t: Q_1 \to Q_0$ called the \textit{source} $s$ and \textit{target} $t$. For every arrow $\alpha \in Q_1$ the maps associate the corresponding source $s(\alpha)$ and target vertex $t(\alpha)$ in $Q_0$. Usually, the quiver is just denoted by $Q$. 

Further using the notation of \cite{MR2197389}, a \textit{path of length} $\ell \geq 1$ with source $a$ and target $b$ is a sequence
\[
	(a \ | \ \alpha_1, \alpha_2, \ldots,\alpha_\ell \ | \ b),
\]
where $\alpha_k \in Q_1$ for all $1 \leq k \leq \ell$ as well as $s(\alpha_1) = a, \ t(\alpha_k) = s(\alpha_{k+1})$ for each $1 \leq k < \ell$ and $t(\alpha_\ell) = b$. Additionally, each vertex $e_i \in Q_0$ has a path of length $\ell = 0$ associated to it, called the \textit{trivial path} at $e_i$ also denoted by $e_i = (e_i \ | \ \ | \ e_i)$.

The path algebra $kQ$ is defined in the usual way with its natural grading 
\[
	kQ = kQ_0 \oplus kQ_1 \oplus kQ_2 \oplus \ldots \oplus kQ_{\ell} \oplus \ldots
\]
where $kQ_\ell$ denotes the $k$-subspace of $kQ$ generated by the set $Q_\ell$ of all paths of length $\ell$.
\end{de}

\begin{de}
For a quiver $Q$, its \textit{double} $\bar{Q}$ is obtained by keeping the vertex set $Q_0$ and adding a new arrow $\alpha^\ast$ from $e_j$ to $e_i$ for each arrow $\alpha \in Q_1$ from $e_i$ to $e_j$.
\end{de}

\begin{de}\label{def23}
Using \cite[Definition 3.1.1.]{MR2335985}, we can define the preprojective algebra associated to a quiver $Q$. Let $Q$ be a quiver and $\bar{Q}$ its double. The \textit{preprojective algebra} is defined as
\[
	A = \Pi_{Q}(k) = k\bar{Q}/\left(\sum_{\alpha \in Q_1}{\alpha \alpha^\ast} - \sum_{\alpha \in Q_1}{\alpha^\ast \alpha}\right).
\]
\end{de}

Before giving the necessary background, let's fix our conventions:
\begin{de}\label{def24}
Let $R = \bigoplus\limits_{j \in \ZZ}{R_j}$ be a graded algebra and $g$ be an automorphism of $R$. 
\begin{itemize}[topsep=0pt]
	\item	We call $g$ a \textit{graded automorphism} if $g(R_i) = R_i$ for all graded pieces $R_i$ of $R$. The \textit{group of all graded automorphisms of $R$} is denoted by $\Aut_{\gr}(R)$.
	\item	A right $R$-module $M$ is called \textit{graded} if $M = \bigoplus_{j \in \ZZ}{M_j}$ and all $r \in R_n$ and $m \in M_\ell$ satisfy $m r \in M_{n+\ell}$ for all $n, \ell \in \ZZ$. 
	\item	For a graded $R$-module $M = \bigoplus_{j \in \ZZ}{M_j}$ a linear map $\varphi$ satisfying $\varphi(M_j) \subseteq M_j$ is called \textit{graded}. The set of all graded automorphisms of $M$ is denoted by $\Aut_{\gr}(M)$.
	\item	For a graded $R$-module $M = \bigoplus_{j \in \ZZ}{M_j}$, the \textit{shifted module} $M[\ell] = \bigoplus_{j \in \ZZ}{(M[\ell])_j}$ is defined by $(M[\ell])_j = M_{\ell+j}$.
	\item	For a right $R$-module $M$, the \textit{$g$-twisted module} $M^g$ is defined as $M^g = M$ with changed action $m \ast r = m g(r)$ for $m \in M, \ r \in R$.
	\item	Let $M$ and $N$ be graded right $R$-modules. A $k$-linear map $\varphi: M \to N$ is called \textit{$g$-linear} if $\varphi(m r) = \varphi(m) g(r)$ for all $m \in M$ and $r \in R$. As noted in \cite{MR1438180}, a graded map $\varphi: M \to N$ is $g$-linear if and only if it is a graded $R$-homomorphism from $M$ to $N^g$.
\end{itemize}
\end{de}

Inspired by \cite[Lemma 1.5.5]{DanNotes}, this following standard lemma becomes handy later when we compute $\Ext^j((e_iA)/(e_iA)_{\geq 1},A)$. We omit the straightforward proof.

\begin{lemma}\label{lemma25}
Let $B$ be a noetherian $\NN$-graded algebra with degree zero piece $B_0 \cong k^n$. Denote the primitive orthogonal idempotents by $\{e_1,\ldots,e_n\}$. Also, denote the indecomposable projective right $B$-modules by $P_i = e_iB$. Recall that $\Hom(e_iB, e_jB) \cong e_jBe_i$ via left multiplication.
\begin{itemize}[topsep=0pt]
	\item[(a)]	For $s_j \in \{1,\ldots,n\}$ for $j = 1,\ldots,N$, there exists a left $B$-module isomorphism
	\[
		\Hom_B\left(\bigoplus_{j = 1}^{N}{P_{s_j}}, B\right) \cong \bigoplus_{j=1}^{N}{B e_{s_j}}.
	\]
	\item[(b)]	Given a graded right $B$-module homomorphisms $\Gamma: P = \bigoplus_{i = 1}^{N}{P_{s_i}} \to Q = \bigoplus_{j = 1}^{M}{P_{t_j}}$, represented by left multiplication by a matrix $T$, i.e. $\Gamma(\beta) = T \cdot \beta$, whose entries $T_{ji}$ lie in $e_{t_j}Be_{s_i}$. Applying the functor $\Hom(-,B)$ yields a left module homomorphism $\Gamma^\ast: \Hom(Q,B) \to \Hom(P,B)$. By (a), 
	\[
		\Gamma^\ast: \bigoplus_{j=1}^{M}{B e_{t_j}} \to \bigoplus_{i = 1}^{N}{B e_{s_i}}
	\]
	and $\Gamma^\ast$ is given by right multiplication by $T$.
\end{itemize}
\end{lemma} 

\begin{de}
Let $R$ be a noetherian $\NN$-graded algebra. A graded right $R$-module $M = \bigoplus_{j \in \ZZ}{M_j}$ is called \textit{bounded} if there exist $c_1$, $c_2 \in \ZZ$ such that $M_j = 0$ for all $j < c_1$ and $M_j = 0$ for all $j > c_2$. 
\end{de}

A weak homological property introduced by Artin and Zhang is the $\chi$-condition for a noetherian $\NN$-graded algebra from \cite[Definitions 3.2]{MR1304753}.
\begin{de}\label{defchi}
A noetherian $\NN$-graded algebra $R$ satisfies the \textit{$\chi$-condition} if $\Ext_R^j(R_0,M)$ is bounded for all $j$ and all finitely generated right $R$-modules $M$. 
\end{de}

Our definition of the $\chi$-condition is not the original definition given in \cite[Definition 3.7]{MR1304753}. However, due to \cite[Proposition 3.11 (2)]{MR1304753} the original definition is equivalent to Definition \ref{defchi} as long as $R$ is locally finite which we assume for all our algebras.

\begin{lemma}\label{lemma27}
Let $B$ be a noetherian $\NN$-graded algebra with $B_0 \cong k^n$. Denote the primitive orthogonal idempotents by $\{e_1,\ldots,e_n\}$. If $B$ has finite global dimension $d$ and $\Ext^i_B(S_j,B) \cong \delta_{di} S_{\sigma(j)}^\ast$ for some $\sigma \in \Sym(n)$ and $S_j = e_j B / (e_j B)_{\geq 1}, \ S_j^\ast = Be_j/(Be_j)_{\geq 1}$ for all $j$, then $B$ satisfies the $\chi$-condition.
\begin{pro}
Let $M$ be a finitely generated graded $B$-module. It suffices to show that $\Ext^i(B_0,M)$ is bounded for all $i$. We proceed by induction on $\pdim(M)$. If $\pdim(M) = 0$, then $M$ is projective and therefore $M$ is isomorphic to a finite sum of shifted copies of $P_\ell = e_\ell B$ for $\ell = 1,\ldots,n$. Using $B_0 \cong \bigoplus{S_j}$ and the hypothesis, $\Ext^i(B_0,B)$ is bounded. Therefore, its summands $\Ext^i(B_0,e_\ell B)$ are bounded for $\ell = 1,\ldots, n$. Being a sum of shifted copies of the $e_\ell B$'s, also for $M$ it is true that $\Ext^i(B_0,M)$ is bounded. If $\pdim(M) > 0$, there is a short exact sequence
\[
	0 \longrightarrow N \longrightarrow P \longrightarrow M \longrightarrow 0
\]
where $P$ is projective and $P \longrightarrow M$ is a minimal surjection. This means $P \longrightarrow M$ is the start of a minimal projective resolution of $M$ and it follows that $\pdim(N) = \pdim(M) - 1$. By induction, $\Ext^i(B_0,N)$ and $\Ext^i(B_0,P)$ are bounded for all $i$. Consider the long exact sequence
\[
\begin{tikzcd}[row sep = 0ex
    ,/tikz/column 1/.append style={anchor=base west}
    ,/tikz/column 2/.append style={anchor=base west}
    ,/tikz/column 3/.append style={anchor=base west}
    ,/tikz/column 4/.append style={anchor=base west}
    ]
\ldots \arrow[r] & \Ext^i(B_0,N) \arrow[r] & \Ext^i(B_0,P) \arrow[r] & \Ext^i(B_0,M) \\
\ \ \arrow[r] & \Ext^{i+1}(B_0,N) \arrow[r] & \Ext^{i+1}(B_0,P) \arrow[r] & \Ext^{i+1}(B_0,M) \arrow[r] & \ldots \\
\ \ \arrow[r] & \Ext^d(B_0,N) \arrow[r] & \Ext^d(B_0,P) \arrow[r] & \Ext^d(B_0,M) \arrow[r] & 0,
\end{tikzcd}
\]
where $d = \gldim(B) < \infty$. This shows that $\Ext^i(B_0,M)$ is always squeezed between two bounded $\Ext$'s and therefore, $\Ext^i(B_0,M)$ is also bounded for all $i$.\qed
\end{pro}
\end{lemma}

Finally, let us introduce the term twisted Calabi-Yau. The proposition afterwards will be used to show that the fixed ring is twisted Calabi-Yau by only proving the generalized Gorenstein condition together with having finite global dimension. Moreover, it helps completing the list of properties of $A$ in Proposition \ref{prop210}.
\begin{de}[{\cite[Definition 1.2]{DanToDo}}]\label{CYdef}
Let $B$ be a $k$-algebra and call $B^e = B \otimes_k B^{\op}$. We say that $B$ is \textit{twisted Calabi-Yau of dimension $d$} if
\begin{itemize}[topsep=0pt, partopsep=0pt,leftmargin=0.75in]
	\item[(CY 1)] $B$ has a resolution of finite length by finitely generated projective $(B,B)$-bimodules and
	\item[(CY 2)] there exists an invertible $k$-central $(B,B)$-bimodule $U$ such that
	\[
		\Ext_{B^e}^i(B,B^e) \cong \begin{cases} 0, & i \neq d \\ U, & i = d \end{cases}
	\]
	as $(B,B)$-bimodules, where each $\Ext_{B^e}^i(B,B^e)$ is considered as a right $B^e$-module via the ``inner" right $B^e$-structure of $B^e$.
\end{itemize}
\end{de}

\begin{prop}[{\cite[Theorem 5.15]{DanToDo}}]\label{prop29}
Let $B$ be an $\NN$-graded algebra with $B_0 \cong k^n$. Denote the primitive orthogonal idempotents by $\{e_1,\ldots,e_n\}$, the simple right $B$-modules by $S_j = e_jB/(e_jB)_{\geq 1}$ and the simple left $B$-modules by $S_j^\ast = Be_j/(Be_j)_{\geq 1}$, respectively. Then the following are equivalent:
\begin{itemize}[topsep=0pt]
	\item[(1)]	$B$ is twisted Calabi-Yau of dimension $d$.
	\item[(2)]	$\gldim(B) = d$ and $B$ satisfies $\Ext_B^i(S_j,B) \cong \delta_{di}S_{\sigma(j)}^\ast[\ell_j]$ for some $\sigma \in \Sym(n)$ and some $\ell_j \in \NN$ for all $j$. This is called the generalized Gorenstein condition.
\end{itemize}
\end{prop}

\begin{prop}\label{prop210}	Here are some facts about a preprojective algebra $A = \Pi_Q(k)$ of an extended Dynkin quiver $Q$ with vertex set $\{e_1,\ldots,e_n\}$. 
 \begin{itemize}[topsep=0pt]
	\item[(a)]	$A$ has global dimension $2$ and satisfies the generalized Gorenstein condition 
	\[
		\Ext^i_A(S_j,A) \cong \delta_{di} S_{\sigma(j)}^\ast[\ell_j]
	\]
	for the simple right modules $S_j = e_jA/(e_jA)_{\geq 1}$ for all $j = 1,\ldots,n$ and $\sigma = \id$. 
	\item[(b)]  $A$ is Calabi-Yau of dimension $d = 2$.
	\item[(c)]	Lemma \ref{lemma27}. $A$ satisfies the $\chi$-condition from Definition \ref{defchi}.
	\item[(d)]	\cite[Theorem]{MR624903}. $\GK(\Pi_Q(k)) = 2$. (The typo was noted in \cite[Proposition 6.10]{MR876985}).
	\item[(e)]	\cite[Corollary 5.17]{MR2197389}. The trivial paths $\{e_1,\ldots,e_n\}$ associated to the vertices form the set of primitive and pairwise orthogonal idempotents. Then, the $P_i = e_iA$ for $i = 1,\ldots,n$ are the only indecomposable projectives and thus every finitely generated projective module is a sum of finitely many possibly shifted copies of the $P_i$'s.
	\item[(f)]	\cite[Theorem 6.5]{MR876985}. $A$ is a noetherian polynomial identity ring.
\end{itemize}
\begin{pro}
By Proposition \ref{prop29} it follows that showing (a) also proves (b). In other words, we need to show that $d = \gldim(A) = 2$ and that $\Ext^i_A(S_j,A) \cong \delta_{di} S_j^\ast = \delta_{di}Ae_j/(Ae_j)_{\geq 1}$ for all simple modules $S_j = e_jA/(e_jA)_{\geq 1}$ for $j = 1,\ldots,n$. This requires minimal projective resolutions of the $S_j$'s.

Consider $Q = \widetilde{A_{n-1}}$ first - for basic information on $A = \Pi_Q(k)$ look at Remark (3.1). The proof of \cite[Theorem 3.2]{MR2355031} gives the first three steps of the minimal projective resolution $(P_j)_\bullet$ of $S_j$ which is also where it stops:
\[
	0 \longrightarrow e_jA \overset{\delta_2 = \begin{pmatrix} \alpha_j^\ast \\ -\alpha_{j-1}\end{pmatrix} \boldsymbol{\cdot}}{\xrightarrow{\hspace*{3cm}}} e_{j+1}A \oplus e_{j-1}A \overset{\delta_1 = \begin{pmatrix} \alpha_j & \alpha_{j-1}^{\ast} \end{pmatrix} \boldsymbol{\cdot}}{\xrightarrow{\hspace*{3cm}}} e_jA \overset{\delta_0}{\longrightarrow} S_j \longrightarrow 0,
\]
where $j = 1,\ldots,n$ and we understand all indices mod $n$. This shows that $\gldim(A) = 2 = d$ by \cite[Proposition 3.19 (3)]{DanToDo}. Applying $\Hom(-,A)$ to $(P_j)_\bullet$, Lemma \ref{lemma25} yields
\[
	0 \longrightarrow Ae_j \overset{\delta_1^\ast = \ \boldsymbol{\cdot} \begin{pmatrix} \alpha_j & \alpha_{j-1}^\ast\end{pmatrix}}{\xrightarrow{\hspace*{3cm}}} Ae_{j+1} \oplus Ae_{j-1} \overset{\delta_2^\ast = \ \boldsymbol{\cdot} \begin{pmatrix} \alpha_j^\ast \\ -\alpha_{j-1} \end{pmatrix}}{\xrightarrow{\hspace*{3cm}}} Ae_j \longrightarrow 0.
\]
Again from the proof of \cite[Theorem 3.2]{MR2355031} it follows that this is the start of a minimal projective resolution of $S_j^\ast$. Since there is no relation of the form $\beta \alpha_j = 0$ for some nonzero $\beta \in Ae_j$ we know that $\delta_1^\ast$ is injective. These two observations yield $\Ext^2_A(S_j,A) = S_j^\ast$ and $\Ext^i_A(S_j,A) = 0$ for $i \neq 2$ which shows (a). Proposition \ref{prop29} implies (b). The cases $\widetilde{D_{n-1}}$ and $\widetilde{E_m}$ for $m = 6,7,8$ are similar and can be found in the author's dissertation \cite{dissertation}.\qed
\end{pro}
\end{prop}

Now, we can recall the definitions of the Hilbert series and the trace function.

\begin{de}\label{def211}
Let $R$ be an $\NN$-graded algebra over $k$ and let $M$ be a left bounded graded locally finite right $R$-module. Further, let $g \in \Aut_{\gr}(M)$ be a graded automorphism of $M$. Then:
\begin{itemize}[topsep=0pt]
	\item	The \textit{Hilbert series} $\HS_M(t)$ of $M$ is the formal Laurent series
\[
	\HS_M(t) = \sum_{j \in \ZZ}{\dim_k(M_j)t^j}.
\]	
	\item	As in \cite{MR1438180}, the \textit{trace} $\Tr_M(g,t)$ of $g$ is the formal Laurent series
\[
	\Tr_M(g,t) = \sum_{j \in \ZZ}{\tr_k(g|_{M_j})t^j}. 
\]
\end{itemize}
\end{de}

\begin{de}\label{def212}
Let $B$ be a noetherian $\NN$-graded algebra with $B_0 \cong k^n$. Denote the primitive orthogonal idempotents by $\{e_1,\ldots,e_n\}$. Let $g \in \Aut_{\gr}(B)$ be a graded automorphism of $B$ such that $g(e_i) = e_i$ for every $i = 1,\ldots,n$. Notice that $e_iB$ is a projective right $B$-module and for each $e_i$, there is an induced map 
\[
	g_i: e_iB \longrightarrow e_iB, \ \ e_i \beta \mapsto e_i g(\beta),
\]
denoted by $g_i$. Further, $g_i = g|_{e_iB}$ is a $g$-linear right $B$-module isomorphism from $e_iB$ to $e_iB$.
\end{de}

\begin{de} 
Let $B$ be a noetherian $\NN$-graded algebra with $B_0 \cong k^n$. Denote the primitive orthogonal idempotents by $\{e_1,\ldots,e_n\}$. The \textit{matrix Hilbert series} $\HM_{B}(t)$ or $\HM(t)$ is the matrix with entries
\[
	[\HM_{B}(t)]_{ij} = [\HM(t)]_{ij} = \sum_{s = 0}^{\infty}{(\dim_k e_i B_s e_j)t^s}. 
\]
For $g \in \Aut_{\gr}(B)$ with $g(e_i) = e_i$, the \textit{matrix trace} $\TrM_B(g,t)$ or $\TrM(g,t)$ (\textit{vector trace} $\TrV_B(g,t)$ or $\TrV(g,t)$, respectively) is the matrix (vector) with entries
\begin{align*}
	[\TrM_B(g,t)]_{ij} &= [\TrM(g,t)]_{ij} = \sum_{s=0}^{\infty}{\tr(g|_{e_i B_s e_j})t^s}, \\
	[\TrV_B(g,t)]_{j} &= [\TrV(g,t)]_{j} = \sum_{s=0}^{\infty}{\tr(g|_{e_jB_s})t^s}.
\end{align*}
\end{de}

\begin{lemma}\label{lemma214}	Let $R$ be an $\NN$-graded ring and $G$ be a finite subgroup of $\Aut_{\gr}(R)$. The fixed ring equals $R^G = \{x \in R \mid g(x) = x, \text{ for all } g \in G\}$. Then:
\begin{itemize}[topsep=0pt]
	\item[(a)]	\cite[Lemma 5.1]{MR1438180}. $R \cong R^G \oplus C$ as $(R^G,R^G)$-bimodules, where $C$ is the kernel of the map
	\[
		\pi_G: x \mapsto \frac{1}{|G|}\sum_{g\in G}{g(x)}.
	\]
	\item[(b)]	\cite[Corollary 1.12]{MR590245}. If $R$ is right noetherian, then $R^G$ is right noetherian. \\
	\item[(c)]	Molien's Theorem, \cite[Lemma 5.2]{MR1438180}. The Hilbert series of the fixed ring equals: 
	\[
		\HS_{R^G}(t) = \frac{1}{|G|} \sum_{g \in G}{\Tr_R(g,t)}.
	\]
\end{itemize}
\end{lemma}

The following is a straightforward generalization of Molien's Theorem, Lemma \ref{lemma214}(c), to the language of matrix Hilbert series and matrix traces.
\begin{thm}\label{thm215}
Let $B$ be a noetherian $\NN$-graded algebra with degree zero piece $B_0 \cong k^n$. Denote the primitive orthogonal idempotents by $\{e_1,\ldots,e_n\}$. Let $G$ be a finite subgroup of $\Aut_{\gr}(B)$ in which every $g \in G$ satisfies $g(e_i) = e_i$ for all $i$. Then 
\begin{align*}
	\HM_{B^G}(t) &= \frac{1}{|G|} \sum_{g\in G}{\TrM_B(g,t)}, \\
	\HS_{e_jB^G}(t) &= \frac{1}{|G|}\sum_{g \in G}{[\TrV_B(g,t)]_{j}}.
\end{align*}
\end{thm}

We conclude this introduction by stating two important results. First, under weak conditions the ring is always finitely generated as a module over its fixed ring and second, the fixed ring satisfies the $\chi$-condition if the ring satisfies this property.

\begin{prop}{(\cite[Corollary 5.9]{MR590245})}\hspace{1em}\label{prop216}
Let $R$ be a right (left) noetherian ring, $G$ be a finite subgroup of $\Aut(R)$ such that $|G|^{-1} \in R$. Then $R$ is a finitely generated right (left) $R^G$-module.
\end{prop}

\begin{prop}{(\cite[Remark after Proposition 8.7]{MR1304753})}\hspace{1em}\label{prop217}
Let $R$ be a noetherian graded ring and $G$ be a finite subgroup of $\Aut_{\gr}(R)$. If $R$ satisfies the $\chi$-condition then so does $R^G$.
\end{prop}

\begin{other}{Convention}
Unless said otherwise, $A = \Pi_Q(k)$ for an extended Dynkin quiver $Q$ with vertices $\{e_1,\ldots,e_n\}$. All groups $G$ and automorphisms $g$ acting on $A$ have finite order and act by extending scalar multiplication on the arrows of degree $1$. In particular, every $g \in G$ fixes each vertex $e_i$ for $i = 1,\ldots,n$. Also, every automorphism obeys the grading. Lastly, the fixed ring is denoted by $A^G = \{a \in A \mid g(a) = a, \text{ for all } g \in G\}$.
\end{other}

\section{Case $\widetilde{A_{n-1}}$}\label{sectionAn}
Let $n \geq 3$ and let $Q$ be the following quiver coming from the extended Dynkin diagram $\widetilde{A_{n-1}}$:
\[
\begin{tikzcd}
& & & e_{n} \ \bullet \arrow[llldd,"\alpha_{n}"] \\
& \\
e_1 \bullet \arrow[r,swap,"\alpha_1"] & \bullet \ e_2 \arrow[r,swap,"\alpha_2"] & \bullet \ e_3 & \ldots & \bullet \arrow[r] & \bullet \arrow[r] & \arrow[llluu, "\alpha_{n-1}"] \bullet \ e_{n-1}
\end{tikzcd}
\]

\vspace{12pt}
In this section, let $A = \Pi_Q(k)$ be the preprojective algebra with respect to $Q$. We are going to present two formulas in Proposition \ref{prop34} and Proposition \ref{prop35} to compute $\TrV_A(g,t)$ for any graded automorphism $g$ of $A$ which fixes the primitive idempotents. In this situation, we can denote $g(\alpha_i) = c_i \alpha_i$ and $g(\alpha_i^\ast) = t_i \alpha_i^\ast$ and reduce the indices modulo $n$ if needed. Given an automorphism $g$ of finite order, Proposition \ref{prop34} also shows that $\Tr_A(g,t)$ can be written as a rational function where the zeros of the denominator are only roots of unity with $1$ of multiplicity at most $2$. Moreover, the graded automorphism group $\Aut_{\gr}(A)$ is described as a semi-direct product in Proposition \ref{prop36}.

\vspace{12pt}
\begin{other}{Definition \& Remark}\label{remark31}	In Definition \ref{def23}, the preprojective algebra was defined as the path algebra modulo the ideal $I$ generated by the one relation
\begin{align}\label{eqrelA}
\sum_{\alpha \in Q_1}{\alpha \alpha^\ast} - \sum_{\alpha \in Q_1}{\alpha^\ast \alpha}.
\end{align}
By successively multiplying the trivial paths $e_i$ for $i = 1,\ldots,n$ to Equation \eqref{eqrelA}, relations of the form $\alpha_i \alpha_i^\ast = \alpha_{i-1}^\ast \alpha_{i-1}$ for all $i = 1,\ldots,n$ are obtained. Starting now, we think of $I$ being generated by the relations of the form $\alpha_i \alpha_i^\ast = \alpha_{i-1}^\ast \alpha_{i-1}$. As a consequence, we can always rearrange a path in $A$ to first have the nonstar arrows occur followed by the star arrows. The absence of overlap or inclusion ambiguities with respect to the natural ordering $\alpha_1 < \alpha_2 < \dots < \alpha_n < \alpha_1^\ast < \ldots < \alpha_n^\ast$ together with the Diamond Lemma \cite[Theorem 1.15]{DanNotes} says that the images of the reduced words form a $k$-basis. However, the leading terms of the relations are of the form $\alpha_i^\ast \alpha_i$ for all $i$ and so the paths of length $j$ starting at $e_\ell$ first having the nonstar arrows occur followed by the star arrows form a basis of $(e_\ell A)_j$. This implies that $\dim((e_\ell A)_j) = j+1$ and therefore
\[
	\HS_{e_\ell A}(t) = 1 + 2t + 3t^2 + 4t^3 + \ldots = \sum_{j = 0}^{\infty}{(j+1)t^j} = \frac{1}{(1-t)^2}. 
\]
Since $A = \bigoplus\limits_{\ell = 1}^{n}{e_\ell A}$, this shows that $\HS_A(t) = \dfrac{n}{(1-t)^2}$.

Moreover, the relations $\alpha_i \alpha_i^\ast = \alpha_{i-1}^\ast \alpha_{i-1}$ force $c_it_i = c_{i-1}t_{i-1}$ for all $i = 1,\ldots,n$. In other words, $c_1t_1 = c_it_i$ for all $i$ is necessary to obtain an automorphism $g$. On the other hand, for every choice of nonzero scalars $c_1,\ldots,c_n, \ t_1,\ldots,t_n$ with $c_1t_1 = c_it_i$ one defines a graded automorphism of $A$.

Every element of $A$ can be represented as a $k$-linear combination of paths in the double of $Q$. In order to simplify notation, we identify an element of $A$ with any representative in $k\bar{Q}$. We call an element of the form $\prod_{j}{\alpha_{i_j}}\prod_{v}{\alpha_{u_v}^\ast}$ a \textit{simple path}. Due to the relations, every path has a unique representation as a simple path. Moreover, every simple path in $A$ is uniquely determined by its starting vertex, its length and the number of star arrows occurring. Further, every element of $A$ is a unique linear combination of simple paths.
\end{other}

\begin{de}
We call a path \textit{purely nonstar}, if it consists of only nonstar arrows and \textit{purely star} if it consists of only star arrows. We call a path \textit{pure} if it is either purely nonstar or purely star. We call a path \textit{mixed} if it is not pure.
\end{de}

\begin{de}
The \textit{type} $\type(\beta) = (s,t)$ of a path $\beta \in A$ is defined by $s = \#$ of nonstar arrows in $\beta$ and $t = \#$ of star arrows in $\beta$. Thanks to the relations in $A$ this definition is well-defined. Two paths $\beta_1$ and $\beta_2$ are \textit{of the same type} if $\type(\beta_1) = \type(\beta_2)$.
\end{de}
Now, we are ready to prove the main proposition of this section which gives valuable information about the vector trace and the trace of a graded automorphism of $A$.

\begin{prop}\label{prop34} Let $A = \Pi_Q(k)$ and let $g$ be a graded automorphism of $A$ with $g(e_i) = e_i$ for all $i = 1,\ldots,n$ of not necessarily finite order. Therefore, we can denote $g(\alpha_i) = c_i \alpha_i$ and $g(\alpha_i^\ast) = t_i \alpha_i^\ast$. Also, recall $g_i = g|_{e_iA}$ and fix the convention that an empty product equals $1$. Then we have the following:
\begin{itemize}[topsep=0pt]
	\item[(a)]	The vector $(1-t_1 \cdots t_n t^n) \cdot \TrV_A(g,t)$ equals	
	\begin{align*}
		\begin{pmatrix} 1 & -c_1 t & & & & & \\ & 1 & -c_2 t & & & & \\ & & 1 & -c_3 t & & & \\ & & & \ddots & \ddots & & \\ & & & & \ddots & \ddots & \\ & & & & & 1 & -c_{n-1} t \\ -c_n t & & & & & & 1 \end{pmatrix}^{-1} \cdot \begin{pmatrix} \sum_{k = 0}^{n-1}{\left(\prod_{j = n-k+1}^{n}{t_j}\right)}t^k \\[1ex] \sum_{k = 0}^{n-1}{\left(\prod_{j = n-k+2}^{n+1}{t_j}\right)}t^k \\[1ex] \sum_{k = 0}^{n-1}{\left(\prod_{j = n-k+3}^{n+2}{t_j}\right)}t^k \\[1ex] \sum_{k = 0}^{n-1}{\left(\prod_{j = n-k+4}^{n+3}{t_j}\right)}t^k \\[1ex] \vdots \\[1ex] \sum_{k = 0}^{n-1}{\left(\prod_{j = 2n-k}^{2n-1}{t_j}\right)}t^k \end{pmatrix}.
\end{align*}
	\item[(b)]	The trace $\Tr_A(g,t)$ can be written as fraction $\frac{p(t)}{q(t)}$ of the two polynomials
\begin{align*}
	p(t) &= \sum_{v=0}^{2n-2}{\left(\sum_{\substack{v = k+s \\ 0 \leq k,s \leq n-1}}{\sum_{\ell = 1}^{n}{\prod_{j=n-k+\ell}^{n+\ell-1}{t_j}\prod_{r = n-s+\ell}^{n+\ell-1}{c_r}}}\right)t^v}, \\
	q(t) &= (1-c_1 \cdots c_n t^n)(1-t_1 \cdots t_n t^n).
\end{align*}
Notice that $p(t)$ and $q(t)$ may not be coprime.
	\item[(c)]	Similarly, $\Tr_{e_iA}(g_i,t) = \frac{p_i(t)}{q_i(t)}$ with $q_i(t) = (1-c_1 \cdots c_n t^n)(1-t_1 \cdots t_n t^n)$.
	\item[(d)]	If $g$ has finite order, all roots of $q(t)$ ($q_i(t)$, respectively) are roots of unity. Additionally, the order of $1$ as a pole of $\Tr_A(g,t)$ ($\Tr_{e_iA}(g_i,t)$, respectively) can be at most two.
\end{itemize}
\begin{pro}
First, recall $g_i: e_iA \to e_iA$, $g_i(e_i \alpha) = g(\alpha)$ from Definition \ref{def212}. Since the primitive idempotents are fixed, this is a well-defined vector space automorphism of $e_iA$. Understanding $\alpha_i$ as a path in degree $1$ from $e_i$ to $e_{i+1}$ (indices are taken modulo $n$), the quotient vector space $B_i = e_iA/\alpha_ie_{i+1}A$ arises naturally. Further using that $g_i$ acts diagonally on $\alpha_i$, the induced maps
\[
	\tilde{g_i}: B_i \to B_i, \ \tilde{g_i}(\beta + \alpha_i e_{i+1}A) = g_i(\beta)+ \alpha_i e_{i+1}A \in B_i
\]
make sense. As noticed in the remark from the beginning of this section, the Diamond Lemma provides a basis of $e_{i+1}A$ consisting of simple paths. This form of the basis and the relations guarantee that the map $\gamma \mapsto \alpha_i \gamma$ for $\gamma \in e_{i+1}A$ is injective. Therefore, due to the trace only seeing the vector space structure, we obtain
\begin{align}\label{eq2}
	\Tr_{B_i}(\tilde{g_i},t) = \Tr_{e_iA}(g_i,t) - c_it\Tr_{e_{i+1}A}(g_{i+1},t).
\end{align}
To proceed, $\Tr_{B_i}(\tilde{g_i},t)$ has to be understood better. Every path in $e_iA$ is a linear combination of simple paths starting at $e_i$. Due to the relations in $A$, a simple path is uniquely determined by its starting vertex and its type. Further, every simple path $\beta \in e_iA$ which has at least one nonstar arrow can be rewritten as $\beta = \alpha_i \beta_1 \in \alpha_ie_{i+1}A$. Therefore, the only nonzero elements of $B_i$ are idempotents and products of star arrows. This provides the following $k$-basis for $B_\ell$:
\begin{align*}
	&\overline{e_\ell}, \ \overline{\alpha_{\ell-1}^\ast}, \ \overline{\alpha_{\ell-1}^\ast\alpha_{\ell-2}^\ast}, \ \overline{\alpha_{\ell-1}^\ast\alpha_{\ell-2}^\ast \alpha_{\ell-3}^\ast}, \ \overline{\alpha_{\ell-1}^\ast\alpha_{\ell-2}^\ast \alpha_{\ell-3}^\ast \alpha_{\ell-4}^\ast}, \ \ldots, \overline{\alpha_{\ell-1}^\ast \alpha_{\ell-2}^\ast \cdots \alpha_1^\ast \alpha_n^\ast \alpha_{n-1}^\ast \cdots \alpha_{\ell+1}^\ast}\\
	&\overline{\alpha_{\ell-1}^\ast \alpha_{\ell-2}^\ast \cdots \alpha_1^\ast \alpha_n^\ast \alpha_{n-1}^\ast \cdots \alpha_{\ell+1}^\ast \alpha_{\ell}^\ast}, \ \overline{\alpha_{\ell-1}^\ast \alpha_{\ell-2}^\ast \cdots \alpha_1^\ast \alpha_n^\ast \alpha_{n-1}^\ast \cdots \alpha_{\ell+1}^\ast \alpha_{\ell}^\ast \alpha_{\ell-1}^\ast}, \ldots
\end{align*}
On these basis elements, $\tilde{g_\ell}$ acts via scalar multiplication by $1, \ t_{\ell-1}, \ t_{\ell-1}t_{\ell-2}, \ t_{\ell-1}t_{\ell-2}t_{\ell-3},$ and so on. Notice that $\tilde{g_\ell}$ applied to $\overline{\alpha_{\ell-1}^\ast \alpha_{\ell-2}^\ast \cdots \alpha_1^\ast \alpha_n^\ast \alpha_{n-1}^\ast \cdots \alpha_{\ell+1}^\ast \alpha_{\ell}^\ast \alpha_{\ell-1}^\ast}$ causes a scalar multiplication by $t_1 \cdots t_n t_{\ell-1}$. This repetition continues and so the trace computes as
\begin{equation}
\begin{aligned}\label{eq3}
	\Tr_{B_\ell}(\tilde{g_\ell},t) &= \frac{1 + t_{\ell-1}t + t_{\ell-1}t_{\ell-2}t^2 + t_{\ell-1}t_{\ell-2}t_{\ell-3}t^3 + \ldots + t_{\ell - 1} t_{\ell-2} \cdots t_1 t_n t_{n-1} \cdots t_{\ell + 1}t^{n-1}}{(1 - t_1 \cdots t_n t^n)} \\ 
	&= \frac{\sum_{k = 0}^{n-1}{\left(\prod_{j = n-k+\ell}^{n+\ell-1}{t_j}\right)}t^k}{(1-t_1 \cdots t_n t^n)}.
\end{aligned}
\end{equation}
Now, we are ready to apply matrix notation to the formulas from equation \eqref{eq2}:
\begin{align}\label{eq4}
	\begin{pmatrix}
		\Tr_{B_1}(\tilde{g_1},t) \\[1ex] \Tr_{B_2}(\tilde{g_2},t) \\[1ex] \Tr_{B_3}(\tilde{g_3},t) \\[1ex] \Tr_{B_4}(\tilde{g_4},t) \\[1ex] \vdots \\[1ex] \Tr_{B_n}(\tilde{g_n},t) \end{pmatrix} &= \underbrace{\begin{pmatrix} 1 & -c_1 t & & & & & \\ & 1 & -c_2 t & & & & \\ & & 1 & -c_3 t & & & \\ & & & \ddots & \ddots & & \\ & & & & \ddots & \ddots & \\ & & & & & 1 & -c_{n-1} t \\ -c_nt & & & & & & 1 \end{pmatrix}}_{=: C^{-1}} \cdot \begin{pmatrix} \Tr_{e_1A}(g_1,t) \\[1ex] \Tr_{e_2A}(g_2,t) \\[1ex] \Tr_{e_3A}(g_3,t) \\[1ex] \Tr_{e_4A}(g_4,t) \\[1ex] \vdots \\[1ex] \Tr_{e_nA}(g_n,t) \end{pmatrix}.
\end{align}
Combining the equations \eqref{eq3} and \eqref{eq4} finishes the proof of part (a).

For part (b), Cramer's rule helps inverting the matrix $C^{-1}$. It says that the entries of the matrix $C = (c_{ij})_{i,j}$ equal $\left(\frac{1}{\det(C^{-1})} \adj(C^{-1})\right)_{i,j}$ where $\adj(-)$ denotes the adjugate matrix. Since $C^{-1}$ is very sparse, $\det(C^{-1})$ can be easily calculated by expanding along the first row. After the first step, one faces a constant term $1$ plus $(-1)^{n} c_1t \cdot (-c_nt)$ multiplied by the determinant of the lower triangular matrix with diagonal entries $(-c_1t,\ldots,-c_{n-1}t)$. This gives $(1-c_1\cdots c_nt^n)$. Now, the $(i,j)$-th entry of $\adj(C^{-1})$ equals $(-1)^{i+j}$ times the determinant of the minor of $C^{-1}$ obtained by deleting the $j$-th row and the $i$-th column. The diagonal entries of $\adj(C^{-1})$ are easily seen to equal $1$ since deleting an $i$-th row and column always isolates a diagonal entry $1$ to be the only entry in a particular column. Every step produces a new column with just one entry on the diagonal equal to $1$. As for $i < j$, expanding along the columns bigger than $j$ and the rows smaller than $i$ which all only have one diagonal entry equal to $1$ reduces the problem to a lower triangular matrix whose diagonal entries are $(-c_it,\ldots,-c_{j-1}t)$. A similar analysis for $i > j$ finally gives
\begin{align*}
	(1 - c_1 \cdots c_n t^n) \cdot c_{ii} &= 1, \\
	(1 - c_1 \cdots c_n t^n) \cdot c_{ij} &= \prod_{k = i}^{j-1}{(c_k t)},	\hspace{36pt} i < j, \ \ \ \text{ and} \\
	(1 - c_1 \cdots c_n t^n) \cdot c_{ij} &= \prod_{k = i}^{n+j-1}{(c_k t)},	\hspace{25pt} i > j.
\end{align*}
This allows us to multiply $C$ from the left on both sides of Equation \eqref{eq4}. Being interested in $\Tr_A(g,t)$ means to multiply $C$ from the left to the vector $(\Tr_{B_k}(\tilde{g_k},t))_{k = 1,\ldots,n}$ and sum up the entries of the resulting vector. However, the same result can be obtained by first summing up each column of $C$ to obtain a vector $\col(C)$ and then computing the dot product of $\col(C)$ and $(\Tr_{B_k}(\tilde{g_k},t))_{k = 1,\ldots,n}$. The $j$-th entry of $\col(C)$ equals:
\begin{small}
\begin{align}\label{eq5}
	(1 - c_1 \cdots c_n t^n) \cdot \col(C)_j &= \left(1 + \sum_{i = 1}^{j-1}{\prod_{k = i}^{j-1}{(c_k t)}} + \sum_{i = j+1}^{n}{\prod_{k = i}^{n+j-1}{(c_k t)}} \right) = \sum_{s = 0}^{n-1}{\left(\prod_{r = n-s+j}^{n+j-1}{c_r}\right)t^s}.
\end{align}
\end{small}

\noindent Finally, putting together the Equations \eqref{eq3} and \eqref{eq5}, the dot product equals
\begin{align*}
\Tr_A(g,t) &= \frac{1}{(1-c_1 \cdots c_n t^n)(1-t_1 \cdots t_n t^n)} \cdot \sum_{\ell = 1}^{n}{\left[\sum_{k = 0}^{n-1}{\left(\prod_{j = n-k+\ell}^{n+\ell-1}{t_j}\right)t^k} \cdot \sum_{s = 0}^{n-1}{\left(\prod_{r = n-s+\ell}^{n+\ell-1}{c_r}\right)t^s}\right]}\\
&= \frac{1}{(1-c_1 \cdots c_n t^n)(1-t_1 \cdots t_n t^n)} \cdot \sum_{v = 0}^{2n-2}{\left[ \sum_{\ell = 1}^{n}{\sum_{\substack{v = k+s\\ 0 \leq k,s \leq n-1}}{\prod_{j = n-k+\ell}^{n+\ell-1}{t_j} \prod_{r = n-s+\ell}^{n+\ell-1}{c_r}}}\right]t^v} \\
&= \frac{1}{(1-c_1 \cdots c_n t^n)(1-t_1 \cdots t_n t^n)} \cdot \sum_{v = 0}^{2n-2}{\left[\sum_{\substack{v = k+s\\ 0 \leq k,s \leq n-1}}{\sum_{\ell = 1}^{n}{\prod_{j = n-k+\ell}^{n+\ell-1}{t_j} \prod_{r = n-s+\ell}^{n+\ell-1}{c_r}}}\right]t^v}.
\end{align*}

Part (c) deals with $\Tr_{e_iA}(g_i,t)$ which equals the $i$-th entry of the vector $\TrV_A(g,t)$. By the proof of part (a) using Cramer's rule, we know that $\Tr_{e_iA}(g_i,t)$ must be a rational function with denominator of the form $q(t) = (1-c_1\cdots c_nt^n)(1-t_1\cdots t_nt^n)$. This is part (c). Regarding part (d), notice that all of $c_1,\ldots,c_n,t_1,\ldots,t_n$ need to be roots of unity such that $g$ has finite order. Therefore, the roots of $1-c_1\cdots c_nt^n$ are roots of unity. Lastly, the maximal order of $1$ as a root of $q(t)$ is two since all roots of $(1-c_1\cdots c_nt^n)$ are distinct. \qed
\end{pro}
\end{prop}

We can now prove a proposition which generalizes the ideas of \cite[Theorem 2.3]{MR1438180}. It gives a different way to calculate the vector trace. In particular, using the following clear formula together with software like \cite[Mathematica, Version 10.4]{Mathematica} provides the best way to effectively calculate traces for the examples in the next section.

\begin{prop}\label{prop35}
Let $A = \Pi_Q(k)$ and let $g$ be a graded automorphism of $A$ with $g(e_i) = e_i$ for all $i = 1,\ldots,n$ of not necessarily finite order. Therefore, we can denote $g(\alpha_i) = c_i \alpha_i$ and $g(\alpha_i^\ast) = t_i \alpha_i^\ast$. Then the vector trace equals:
\[
	\TrV_A(g,t) = \begin{pmatrix} 1+c_1t_1t^2 & -c_1t & & & & -t_nt \\[1ex] -t_1t & 1 + c_1t_1t^2 & -c_2t & & & \\[1ex] & -t_2t & 1+c_1t_1t^2 & -c_3 t & & \\[1ex] & & \ddots & \ddots & \ddots & \\[1ex] & & & \ddots & \ddots & -c_{n-1} t \\[1ex] -c_nt & & & & -t_{n-1}t & 1+c_1t_1t^2 \end{pmatrix}^{-1} \cdot \begin{pmatrix} 1 \\[1ex] 1 \\[1ex] 1 \\[1ex] \vdots \\[1ex] 1 \\[1ex] 1 \end{pmatrix}.
\]
\begin{proof2}
The idea is to concretely calculate the minimal projective resolutions of the $S_j$ and lift the maps $g_j: e_jA \to (e_jA)^g$. There are only two induced maps as the minimal projective resolution is of length $2$. Then, \cite[Lemma 2.1]{MR1438180} gives that the alternating sum of traces of the connecting maps equal zero. All that is left to do is to put the formulas in matrix notation. A detailed proof is given in the author's dissertation \cite{dissertation}. \qed
\end{proof2}
\end{prop}

This closing proposition is independent from the rest of the section. While we will only be working with automorphisms which fix the vertices, it is nice to know that this subgroup of $\Aut_{\gr}(A)$ is of finite index.
\begin{prop}\label{prop36}
The graded automorphism group of $A = \Pi_Q(k)$ for $Q = \widetilde{A_{n-1}}$ is isomorphic to $D_{2n} \ltimes N$, where $N = \{g \in \Aut_{\gr}(A) \mid g(e_i) = e_i, \text{ for } i = 1,\ldots,n\}$ and $D_{2n}$ denotes the dihedral group of order $2n$.
\begin{pro}
We understand $D_{2n}$ as the group of symmetries of the regular $n$-gon whose vertices are denoted by $e_1,\ldots,e_n$. Let $g \in \Aut_{\gr}(A)$. Since $g$ maps a pair of adjacent vertices to a pair of adjacent vertices, there exists an element $d^{-1} \in D_{2n}$ such that $g(e_i) = d^{-1}(e_i)$ for all $i = 1,\ldots,n$. It follows that $(d \circ g)(e_i) = e_i$ for all $i = 1,\ldots,n$. Hence, $d \circ g \in N$. Further, for $s \in D_{2n}$ and $h \in N$ it is true that
\[
	\left(s \circ h \circ s^{-1}\right)(e_j) = (s \circ h)(e_i) = s(e_i) = e_j
\]
whenever $s(e_i) = e_j$. This shows that $s \circ h \circ s^{-1} \in N$ which proves $N$ is a normal subgroup of $\Aut_{\gr}(A)$.\qed
\end{pro}
\end{prop}

\section{Examples}\label{expchapter}
This section provides examples of automorphisms and fixed rings of $A = \Pi_Q(k)$ for $Q = \widetilde{A_2}$ where $A_2$ is the following quiver:
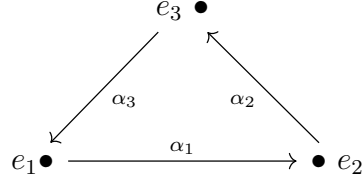
\[
\begin{tikzcd}
& & e_{3} \ \bullet \arrow[ldd,"\alpha_3"] \\
& \\
& e_1 \bullet \arrow[rr,"\alpha_1"] & & \bullet \ e_2 \arrow[luu,"\alpha_2", shift right = 0ex] 
\end{tikzcd}
\]

We use the notation as before, $g \in \Aut_{\gr}(A)$ maps 
\begin{center}
\begin{tabular}{llll}
	$\alpha_1 \mapsto c_1 \alpha_1$, & & & $\alpha_1^\ast \mapsto t_1 \alpha_1^\ast$, \\
	$\alpha_2 \mapsto c_2 \alpha_2$, & & & $\alpha_2^\ast \mapsto t_2 \alpha_2^\ast$, \\
	$\alpha_3 \mapsto c_3 \alpha_3$, & & & $\alpha_3^\ast \mapsto t_3 \alpha_3^\ast$.
\end{tabular}
\end{center}
Remark (3.1) proves $\HS_A(t) = \frac{3}{(1-t)^2}$. We always give an automorphism $g$ and call $G = \langle g \rangle$.

Earlier in Section \ref{sectionAn} was noted that the relations in $A$ are $\alpha_i \alpha_i^\ast = \alpha_{i-1}^\ast \alpha_{i-1}$. Recall that this forces $c_1t_1 = c_2t_2 = c_3t_3$ and whenever $c_1t_1 = c_2t_2 = c_3t_3$ for nonzero $c_i$ and $t_i$, we obtain an automorphism of $A$ by construction.

Proposition \ref{prop34}(b) provides the following formula of the trace of $g$ as
\begin{small}
\begin{align*}
	\Tr_A(g,t) &= \frac{3 + (c_1+c_2+c_3+t_1+t_2+t_3)t + (c_1c_2 + c_1c_3 + c_2c_3 + t_1t_2 + t_1t_3 + t_2t_3 + 3c_1t_1)t^2}{(1-c_1c_2c_3t^3)(1-t_1t_2t_3t^3)} \\
	&\ \ \ + \frac{c_1t_1 (c_1 + c_2 + c_3 + t_1 + t_2 + t_3)t^3 + 3 (c_1t_1)^2t^4}{(1-c_1c_2c_3t^3)(1-t_1t_2t_3 t^3)}.
\end{align*}
\end{small}

Example \ref{exp41} gives a fixed ring of global dimension $2$ coming from an automorphism whose trace has $1$ as a pole of order $1$. A situation that does not happen in the case of Artin-Schelter regular algebras (see \cite[Definition 1.5]{MR2434290} for a definition) manifests in Example \ref{exp42}. Proposition \ref{prop216} shows that $A$ is a finitely generated $A^G$-module and Lemma \ref{lemma214} contributes $A \cong A^G \oplus C$. Thus, for AS-regular algebras, \cite[Lemma 1.10 (c)]{MR2434290} would imply that the fixed ring having finite global dimension implies it is AS-regular. Then \cite[Lemma 1.11 (d)]{MR2434290} says $\gldim(A) = \gldim(A^G)$ and $A$ is free as an $A^G$-module. However, in Example \ref{exp42}, $\gldim(A^G) = 2$ and $A$ is projective but not free as an $A^G$-module. Another situation impossible for AS-regular algebras to occur is revealed in Example \ref{exp43} where the fixed ring has finite global dimension different from the global dimension of the preprojective algebra $\Pi_Q(k)$. 

\begin{ex}\label{exp41}
Consider the following example where $A$ is free over the fixed ring $A^G$. We choose $g$ with $c_1 = 1, \ c_2 = 1, \ c_3 = 1$ and $t_1 = -1, \ t_2 = -1, \ t_3 = -1$. Thanks to the formula in the beginning of the section, the trace of $g$ equals
	\[
		\Tr_A(g,t) = \frac{3}{1-t^2}.
	\]
	Since $c_1 = c_2 = c_3 = 1$, all nonstar arrows are fixed. This implies that a mixed path is fixed if and only if its star part is fixed independently. Hence, it is enough to find the shortest fixed purely star paths. From $t_1 = t_2 = t_3 = -1$ we know that these equal $x = \alpha_3^\ast \alpha_2^\ast, \ y = \alpha_1^\ast \alpha_3^\ast$ and $z = \alpha_2^\ast \alpha_1^\ast$. Using every path's expression as a simple path, we have found all generators of $A^G$. Since a path in $A$ is uniquely determined by its length, its starting vertex and its number of nonstar arrows, the relations 
	\[
		x \alpha_2 = \alpha_1 y, \ \ \ \ \ y\alpha_3 = \alpha_2 z \ \ \ \text{ and } \ \ \ z\alpha_1 = \alpha_3 x
	\]
	appear naturally. Hence, our initial guess is that the fixed ring is the path algebra 
	
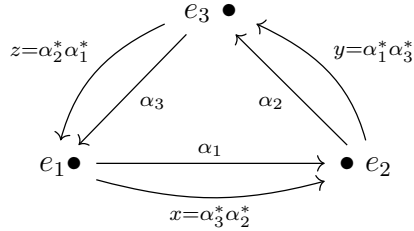
\[
\begin{tikzcd}
& & e_{3} \ \bullet \arrow[ldd,"\alpha_3"] \arrow[ldd,swap,"z = \alpha_2^\ast \alpha_1^\ast", out = 205, in=75, shift right = 0.5ex] \\
& \\
& e_1 \bullet \arrow[rr,swap,"x = \alpha_3^\ast \alpha_2^\ast", out = 345, in=195, shift right = 0.5ex]
 \arrow[rr,"\alpha_1"] & & \bullet \ e_2 \arrow[luu,"\alpha_2", shift right = 0ex] \arrow[luu,swap,"y = \alpha_1^\ast \alpha_3^\ast", out = 105, in=335, shift right = 0.5ex]
\end{tikzcd}
	\] 
	with respect to the stated relations.
	
	The following Gr\"obner basis calculation shows that we actually found a presentation of the fixed ring. We define the order $\alpha_1 < \alpha_2 < \alpha_3 < x < y < z$. Since there are no overlap or inclusion ambiguities, the paths not including the leading terms $x \alpha_2, \ y \alpha_3$ and $z \alpha_1$ form a Gr\"obner basis by the Diamond Lemma. So we need to count them.
	\begin{itemize}[topsep=0pt]
		\item[$\bullet$] Words starting with $x$ can only look like 
			\[
				\begin{tikzcd}[row sep=tiny] 
					x \arrow[r] & x y \arrow[r] & x y z \arrow[r] &	x y z x \arrow[r] & \ldots 
				\end{tikzcd}
			\]
			By symmetry, words starting with $y$ or $z$, respectively, contribute the same to the Hilbert series of $A^G$. Each adds $\dfrac{t^2}{(1-t^2)}$.
		\item[$\bullet$] Words starting with $\alpha_1$ can look like
			\begin{small}
			\[
			\begin{tikzcd}[row sep=tiny]
				\alpha_1 \arrow[dr] \arrow[r] & \alpha_1 y \arrow[r] & \alpha_1 y z \arrow[r] & \alpha_1 y z x \arrow[r] & \alpha_1 y z x y \arrow[r] & \ldots \\
& \alpha_1 \alpha_2 \arrow[dr] \arrow[r] & \alpha_1 \alpha_2 z \arrow[r] & \alpha_1 \alpha_2 z x \arrow[r] &\alpha_1 \alpha_2 z x y \arrow[r] & \ldots \\
& & \alpha_1 \alpha_2 \alpha_3 \arrow[dr] \arrow[r] & \alpha_1 \alpha_2 \alpha_3 x
 \arrow[r] & \alpha_1 \alpha_2 \alpha_3 x y \arrow[r] & \ldots \\
& & & \alpha_1 \alpha_2 \alpha_3 \alpha_1 \arrow[dr] \arrow[r] & \alpha_1 \alpha_2 \alpha_3 \alpha_1 y \arrow[r] & \ldots \\
& & & & \ldots 
\end{tikzcd}
\]\end{small}

\noindent where we see repetition starting at line 4. Combining the first, fourth, seventh and so on row gives $\dfrac{t+t^3+t^5}{(1-t^6)(1-t^3)}$. Similarly, we get $\dfrac{t^2 + t^4 + t^6}{(1-t^6)(1-t^3)}$ and $\dfrac{t^3 + t^5 + t^7}{(1-t^6)(1-t^3)}$ from the second, fifth, eighth and so on line and third, sixth, ninth and so on line, respectively. This adds up to
\[
	\frac{t+t^2+2t^3 + t^4 + 2t^5 + t^6 + t^7}{(1-t^6)(1-t^3)}.
\]
The same has to be added for paths starting with $\alpha_2$ and $\alpha_3$, respectively.
	\item[$\bullet$] Do not forget $e_1, e_2$ and $e_3$ which contribute $3$. 
	\end{itemize}
In the end, we sum up
\begin{align*}
	3 + \dfrac{3t^2}{(1-t^2)} + \frac{3 \cdot (t+t^2+2t^3 + t^4 + 2t^5 + t^6 + t^7)}{(1-t^6)(1-t^3)} = \dfrac{3}{(1-t)^2(1+t)}.
\end{align*}	
	We compare our result to the Hilbert series of $A^G$. According to Molien's Theorem, Lemma \ref{lemma214}(c), $\HS_{A^G}(t)$ equals $\dfrac{3}{(1-t)^2(1+t)}$ which shows that we found a presentation for the fixed ring. It turns out that $A \cong A^G \oplus A^G[-1]$ as right $A^G$-modules via the isomorphism 
	\[
		\varphi: A^G \oplus A^G[-1] \longrightarrow A, \ \ \ \ \ \ (\beta_1,\beta_2) \mapsto \beta_1 + (\alpha_1^\ast + \alpha_2^\ast + \alpha_3^\ast)\beta_2.
	\]
	Linearity is clear. Also, it is easy to see that $\varphi$ is surjective as every simple path with an even number of star arrows is fixed and thus in $A^G$, and every simple path with an odd number of star arrows starting at $e_i$ lives in $\alpha_{i-1}^\ast A^G$. Lastly, $\varphi$ is injective since $\varphi(\beta_1,\beta_2) = 0$ reduces to $\beta_1 = \beta_2 = 0$ as follows. While $\beta_1$ is a linear combination of simple paths with an even number of star arrows $(\alpha_1^\ast + \alpha_2^\ast + \alpha_3^\ast)\beta_2$ can only be composed of simple paths with an odd number of star arrows. Therefore, there is no interaction between the two parts possible and the form of the relations of $A$ do not allow $(\alpha_1^\ast + \alpha_2^\ast + \alpha_3^\ast) \beta_2$ to be zero for a nonzero $\beta_2$.
\end{ex}

In Example \ref{exp41} everything is as for Artin-Schelter regular algebras. The trace of the examined graded automorphism has $1$ as a pole of order $1 = \GK(A) - 1$, the fixed ring $A^G$ is twisted Calabi-Yau and $A$ is a free module over $A^G$. The next examples will show that this does not always happen.

\begin{ex}\label{exp42} This example provides a fixed ring of global dimension $2$ where $A$ is projective but not free over the fixed ring $A^G$. We choose $g$ with $c_1 = 1, \ c_2 = -1, \ c_3 = -1$ and $t_1 = \zeta_3, \ t_2 = -\zeta_3, \ t_3 = -\zeta_3$ where $\zeta_3 = \exp(2 \pi i /3)$. This means $g$ has order $6$ and acts as follows:
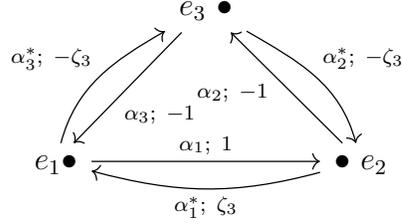
\[
\begin{tikzcd}
& & e_{3} \ \bullet \arrow[ldd,"\alpha_3;\ -1", pos = 0.6] \arrow[rdd,"\alpha_2^\ast; \ -\zeta_3", out=330,in = 105] \\
& \\
& e_1 \bullet \arrow[ruu,"\alpha_3^\ast; \ -\zeta_3", out=75,in = 210] \arrow[rr,"\alpha_1; \ 1"] & & \bullet \ e_2 \arrow[ll,"\alpha_1^\ast; \ \zeta_3", out=195,in = 345] \arrow[luu,"\alpha_2; \ -1", shift right = 0ex, pos = 0.6] 
\end{tikzcd}
	\] 
Thanks to the formula in the beginning of the section and Molien's Theorem, Lemma \ref{lemma214}(c), the trace of the powers of $g$ and the Hilbert series of $A^G$ equal
	\begin{align*}
		\Tr_A(g,t) &= \frac{6 + (-1-i\sqrt{3})t + (-4+4i \sqrt{3}) t^2 + 2t^3 + (-3-3i\sqrt{3})t^4}{2(1-t^3)^2}, \\
		\Tr_A(g^2,t) &= \frac{6+(3-3i\sqrt{3})t}{2(1-t^3)}, \\
		\Tr_A(g^3,t) &= \frac{3-5t+3t^2}{(1-t)(1-t^3)}, \\
		\Tr_A(g^4,t) &= \frac{6+(3+3i\sqrt{3})t}{2(1-t^3)}, \\
		\Tr_A(g^5,t) &= \frac{6 + (-1+i\sqrt{3})t + (-4-4i \sqrt{3}) t^2 + 2t^3 + (-3+3i\sqrt{3})t^4}{2(1-t^3)^2}, \\
		\Tr_A(g^6,t) &= \Tr_A(\id_A,t) = \frac{3}{(1-t)^2}, \ \ \ \ \ \text{ and} \\
		\HS_{A^G}(t) &= \dfrac{3+t+t^2}{(1-t^3)^2}.
	\end{align*}
	Notice that $1$ is a pole of order $1$ of $\Tr_A(g,t)$. In order to find generators for the fixed algebra $A^G$, where $G = \langle g \rangle$, we use the following strategy. Every minimal generator of $A^G$ can have length at most four as every fixed path of length five or larger contains either three nonstar or three star arrows and therefore allows one to decompose the fixed element as concatenation of two individually fixed words. Hence, $\alpha_1, \ x = \alpha_2\alpha_3, \ u_1 = \alpha_3^\ast\alpha_2^\ast \alpha_1^\ast, \ u_2 = \alpha_1^\ast \alpha_3^\ast \alpha_2^\ast, \ s = \alpha_3 \alpha_1 \alpha_2$ and $t = \alpha_2^\ast \alpha_1^\ast \alpha_3^\ast$ are the generators of the fixed ring.

Checking paths of the same type, we find at least the relations $ts = st$, $u_2 x = xu_1$ and $u_1 \alpha_1 = \alpha_1 u_2$ and with respect to the ordering $\alpha_1 < x < u_1 < u_2 < s < t$ there do not exist any overlap or inclusion ambiguities. Therefore, the reduced words form a basis by the Diamond Lemma. To sum up, we know the fixed ring must be isomorphic to a factor of the following path algebra together with the stated relations which we denote by $B$:
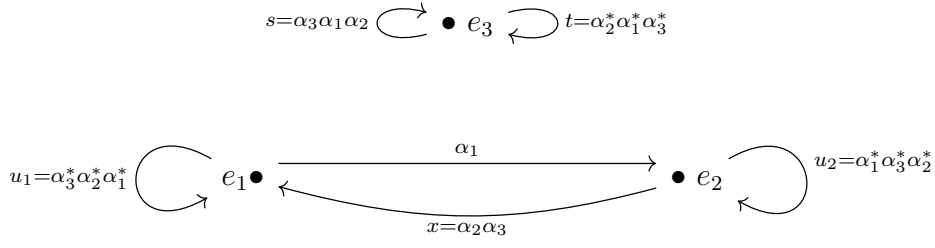
\[
\begin{tikzcd}
& & & \bullet \ e_3 \arrow[loop right, "t = \alpha_2^\ast \alpha_1^\ast \alpha_3^\ast"] \arrow[loop left, "s =\alpha_3 \alpha_1 \alpha_2"] \\
& \\
& e_1 \bullet \arrow[loop left,swap, "u_1 = \alpha_3^\ast \alpha_2^\ast \alpha_1^\ast",out = 150, in = 210, distance=1.5cm] \arrow[rrrr,"\alpha_1", shift left = 1ex]& & & & \bullet \ e_2 \arrow[llll,"x = \alpha_2 \alpha_3", out=195,in = 345, shift left = 0ex] \arrow[loop right, "u_2 = \alpha_1^\ast \alpha_3^\ast \alpha_2^\ast",out = 30, in = 330, distance=1.5cm]
\end{tikzcd}
	\] 
	To verify that this is a presentation of the fixed ring, we count words not containing the leading terms $ts, \ u_2x$ and $u_1 \alpha_1$. 
	\begin{itemize}[topsep = 0pt]
		\item	The part at $e_3$ is isomorphic to a polynomial ring in two variables of degree $3$. Its Hilbert series equals $\HS_{e_3B}(t) = \dfrac{1}{(1-t^3)^2}$.
		\item	The relations imply that the loops $u_i$ for $i = 1,2$ can only be followed up by the same loops again. Including the corresponding idempotents, $\dfrac{1}{(1-t^3)}$ for each $i$ is contributed to the Hilbert series. \\[-8pt]
		\item	Lastly, we need to count paths starting with $\alpha_1$ and $x = \alpha_2 \alpha_3$. Let's start with $\alpha_1$.
		
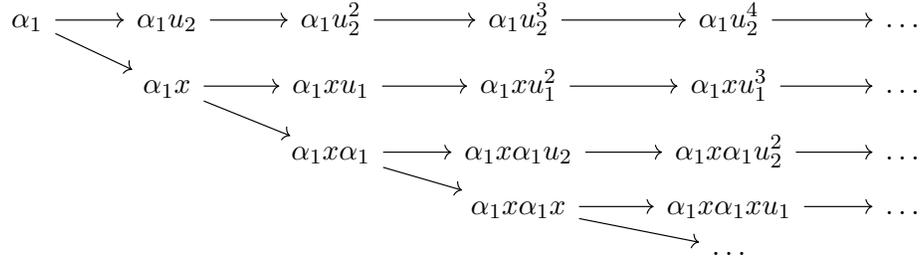
\begin{small}
			\[
			\begin{tikzcd}[row sep=tiny]
				\alpha_1 \arrow[dr] \arrow[r] & \alpha_1 u_2 \arrow[r] & \alpha_1 u_2^2 \arrow[r] & \alpha_1 u_2^3 \arrow[r] & \alpha_1 u_2^4 \arrow[r] & \ldots \\
& \alpha_1 x \arrow[dr] \arrow[r] & \alpha_1 x u_1 \arrow[r] & \alpha_1 x u_1^2 \arrow[r] &\alpha_1 x u_1^3 \arrow[r] & \ldots \\
& & \alpha_1 x \alpha_1 \arrow[dr] \arrow[r] & \alpha_1 x \alpha_1 u_2
 \arrow[r] & \alpha_1 x \alpha_1 u_2^2 \arrow[r] & \ldots \\
& & & \alpha_1 x \alpha_1 x \arrow[dr] \arrow[r] & \alpha_1 x \alpha_1 x u_1 \arrow[r] & \ldots \\
& & & & \ldots 
\end{tikzcd}
\]\end{small}

\noindent Combine all even rows and all odd rows to see the repetition. The odd rows give $\dfrac{t}{(1-t^3)^2}$ and the even rows give $\dfrac{t^3}{(1-t^3)^2}$. For paths starting with $x$ a similar chain of possibilities arises:
			
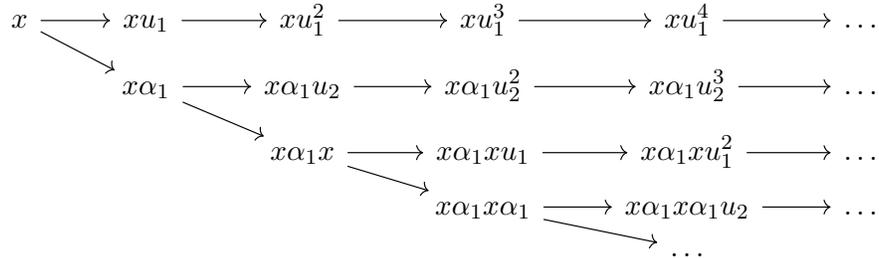
\begin{small}
			\[
			\begin{tikzcd}[row sep=tiny]
				x \arrow[dr] \arrow[r] & xu_1 \arrow[r] & xu_1^2 \arrow[r] & xu_1^3 \arrow[r] & xu_1^4 \arrow[r] & \ldots \\
& x \alpha_1 \arrow[dr] \arrow[r] & x \alpha_1 u_2 \arrow[r] & x \alpha_1 u_2^2 \arrow[r] & x \alpha_1 u_2^3 \arrow[r] & \ldots \\
& & x \alpha_1 x \arrow[dr] \arrow[r] & x \alpha_1 x u_1
 \arrow[r] & x \alpha_1 x u_1^2 \arrow[r] & \ldots \\
& & & x \alpha_1 x \alpha_1 \arrow[dr] \arrow[r] & x \alpha_1 x \alpha_1 u_2 \arrow[r] & \ldots \\
& & & & \ldots 
\end{tikzcd}
\]\end{small}

\noindent Again, combining all the even rows and all the odd rows contributes a total of $\dfrac{t^2+t^3}{(1-t^3)^2}$ to the Hilbert series $\HS_B(t)$.
	\end{itemize}
	Putting the pieces together gives the following sum
	\[
		\HS_B(t) = \frac{1}{(1-t^3)^2} + \frac{2}{(1-t^3)} + \frac{t+t^3}{(1-t^3)^2} + \frac{t^2+t^3}{(1-t^3)^2} = \frac{3+t+t^2}{(1-t^3)^2} = \HS_{A^G}(t)
	\]
	which shows that we have found a presentation for the fixed ring. 	
	
	The previous calculation also helps to find $\HM_{A^G}(t)$. Since the quiver underlying $A^G$ has two connected components, only five entries in the $(3 \times 3)$-matrix are nonzero. The easiest part is $\HS_{e_3A^Ge_3}(t) = \dfrac{1}{(1-t^3)^2}$ as noted above. Hence, the missing entries are $\HS_{e_iA^Ge_j}(t)$ for $i,j = 1,2$. Regarding $\HS_{e_1A^Ge_1}(t)$ (and analogously for $\HS_{e_2A^Ge_2}(t)$) we find as basis $e_1$, powers of $u_1$ and paths starting with $\alpha_1$ and ending in either $x$ or $u_1$. The first two combine to $\dfrac{1}{1-t^3}$. The latter part are exactly the even rows among the paths starting with $\alpha_1$ which contribute $\dfrac{t^3}{(1-t^3)^2}$. This adds up to $\HS_{e_1A^Ge_1} = \dfrac{1}{(1-t^3)^2}$.
	
	In the same way, $\HS_{e_1A^Ge_2}$ has a basis consisting of all the elements starting with $\alpha_1$ and ending in either $u_2$ or $\alpha_1$. These elements are counted by the odd rows starting with $\alpha_1$. From above this means $\HS_{e_1A^Ge_2} = \dfrac{t}{(1-t^3)^2}$.
	
	Analogously, a basis of $e_2A^Ge_1$ consists of the elements of the odd rows starting with $x$. In other words, $\HS_{e_2A^Ge_1} = \dfrac{t^2}{(1-t^3)^2}$. To sum up, we obtain
	\[
		\HM_{A^G}(t) = \frac{1}{(1-t^3)^2} \begin{pmatrix} 1 & t & \\ t^2 & 1 & \\ & & 1 \end{pmatrix}.
	\]
	By \cite[Proposition 3.2.1]{MR2335985}
	\[
		\HM_{A}(t) = (I_2 - Ct + I_2 t^2)^{-1} = \frac{1}{(1-t^3)^2} \begin{pmatrix} 1+t^2+t^4 & t+t^2+t^3 & t+t^2+t^3 \\ t+t^2+t^3 & 1+t^2+t^4 & t+t^2+t^3 \\ t+t^2+t^3 & t+t^2+t^3 & 1+t^2+t^4 \end{pmatrix}
	\]	
	where $C$ is the adjacency matrix of $\bar{Q}$. Inverting the matrix $\HM_{A^G}(t)$, multiplying to $\HM_{A}(t)$ from the right and summing up finally gives
\begin{align*}
	\HS_{e_1A}(t) &= (1+t^2) \HS_{e_1A^G}(t) + t^2 \HS_{e_2A^G}(t) + (t+t^2+t^3)\HS_{e_3A^G}(t), \\
	\HS_{e_2A}(t) &= (t+t^3) \HS_{e_1A^G}(t) + \HS_{e_2A^G}(t) + (t+t^2+t^3)\HS_{e_3A^G}(t), \\
	\HS_{e_3A}(t) &= (t+t^2) \HS_{e_1A^G}(t) + t \HS_{e_2A^G}(t) + (1+t^2+t^4)\HS_{e_3A^G}(t) \ \ \ \ \text{ or } \\[2ex]
	\HS_A(t) &= (1+2t+2t^2+t^3)\HS_{e_1A^G}(t) + (1+t+t^2)\HS_{e_2A^G}(t)\\
	&\ \ \ + (1+2t+3t^2+2t^3+t^4)\HS_{e_3A^G}(t)
\end{align*}
which shows that $A$ is not free as an $A^G$-module. However, Proposition \ref{prop55} will prove that $A$ must be a projective $A^G$-module once we verify $\gldim(A^G) = 2$.

With the notation $P_i = e_iA^G$, the proof of \cite[Theorem 3.2]{MR2355031} states the beginning of the minimal projective resolution of a simple module of a factor of a path algebra. Applied to this example, thanks to computations we know that the minimal projective resolutions of $S_i = P_i/(P_i)_{\geq 1}$ look like
\hspace*{-150pt}
\[
\begin{tikzcd}
0 \arrow[r] & P_2 \arrow[rr, "\begin{pmatrix} u_2 \\ -\alpha_1 \end{pmatrix} \boldsymbol{\cdot}"] & & P_2 \oplus P_1 \arrow[rr, "(\alpha_1 \ \ u_1) \ \boldsymbol{\cdot}"] & & P_1\arrow[r] & S_1 \arrow[r] & 0, \\
0 \arrow[r] & P_1 \arrow[rr, "\begin{pmatrix} u_1 \\ -x \end{pmatrix} \boldsymbol{\cdot}"] & & P_1 \oplus P_2 \arrow[rr, "(x \ \ u_2) \ \boldsymbol{\cdot}"] & & P_2\arrow[r] & S_2 \arrow[r] & 0, & \text{ and} \\
0 \arrow[r] & P_3 \arrow[rr, "\begin{pmatrix} t \\ -s \end{pmatrix} \boldsymbol{\cdot}"] & & P_3 \oplus P_3 \arrow[rr, "(s \ \ t) \ \boldsymbol{\cdot}"] & & P_3\arrow[r] & S_3 \arrow[r] & 0.
\end{tikzcd}
\]
Exactness at the $0$-th and $1$-st position comes for free thanks to \cite[Theorem 3.2]{MR2355031}. At last, each map at the second position is injective because there are no relations of the form $\alpha_1 \beta = 0, \ x \beta = 0$ or $s \beta = 0$ for appropriate nonzero $\beta$'s due to the relations and the chosen ordering. Therefore, \cite[Proposition 3.19 (3)]{DanToDo} provides $\gldim(A^G) = 2$. 
\end{ex}

\begin{ex}\label{exp43} In this example the fixed ring has global dimension $3$. Define $g$ to map
\begin{center}
\begin{tabular}{llll}
	$\alpha_1 \mapsto \zeta_3 \alpha_1$, & & & $\alpha_1^\ast \mapsto \alpha_1^\ast$, \\
	$\alpha_2 \mapsto \alpha_2$, & & & $\alpha_2^\ast \mapsto \zeta_3 \alpha_2^\ast$, \\
	$\alpha_3 \mapsto \zeta_3^2 \alpha_3$, & & & $\alpha_3^\ast \mapsto \zeta_3^2 \alpha_3^\ast,$
\end{tabular}
\end{center}
where $\zeta_3 = \exp(2 \pi i/3)$ and $c_1t_1 = \zeta_3$. This means $g$ acts as follows:
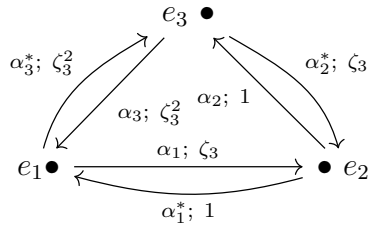
\[
\begin{tikzcd}
& & e_{3} \ \bullet \arrow[ldd,"\alpha_3;\ \zeta_3^2", pos = 0.5] \arrow[rdd,"\alpha_2^\ast; \ \zeta_3", out=330,in = 105] \\
& \\
& e_1 \bullet \arrow[ruu,"\alpha_3^\ast; \ \zeta_3^2", out=75,in = 210] \arrow[rr,"\alpha_1; \ \zeta_3"] & & \bullet \ e_2 \arrow[ll,"\alpha_1^\ast; \ 1", out=195,in = 345] \arrow[luu,"\alpha_2; \ 1", shift right = 0ex, pos = 0.6] 
\end{tikzcd}
	\] 
Thanks to the formula in the beginning of the section and Molien's Theorem, Lemma \ref{lemma214}(c), the traces of $g$ and $g^2$ as well as the Hilbert series of $A^G$ equal
\begin{align*}
\Tr_A(g,t) &= \frac{3i(1-t^2)(-2i+(-i+\sqrt{3})t^2)}{2(1-t^3)^2}, \\
\Tr_A(g^2,t) &= \frac{(-3i)(1-t^2)(2i+(i+\sqrt{3})t^2)}{2(1-t^3)^2}, \ \ \ \ \text{ and} \\
\HS_{A^G}(t) &= \frac{3+2t+2t^2+2t^3}{(1-t^3)^2}.
\end{align*}
The same argument as in the previous example works, i.e. collecting all paths of length smaller than or equal to $4$ which are fixed provides the generators of $A^G$ for $G = \langle g \rangle$. Checking all words of length smaller than or equal to four provides a path algebra $B$ as follows:
	
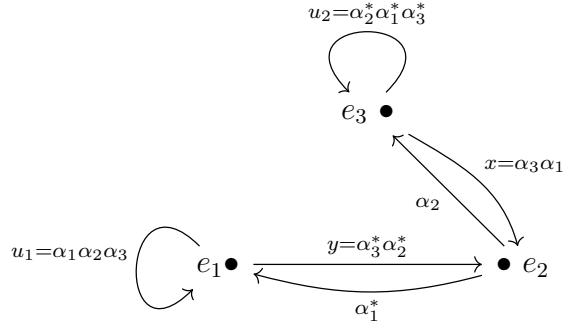
\[
\begin{tikzcd}
& & e_{3} \ \bullet \arrow[loop,swap,"u_2 = \alpha_2^\ast \alpha_1^\ast \alpha_3^\ast"] \arrow[rdd,"x = \alpha_3 \alpha_1", out=330,in = 105] \\
& \\
& \arrow[loop,swap,"u_1 = \alpha_1\alpha_2 \alpha_3",out = 135, in = 225, distance=1.5cm] e_1 \bullet \arrow[rr,"y = \alpha_3^\ast \alpha_2^\ast"] & & \bullet \ e_2 \arrow[ll,"\alpha_1^\ast", out=195,in = 345] \arrow[luu,"\alpha_2", shift right = 0ex] 
\end{tikzcd}
	\]
	
	\vspace{-18pt}
	Notice that $\alpha_1 \alpha_1^\ast \alpha_3^\ast = \alpha_3^\ast \alpha_2^\ast \alpha_2$. We fix the ordering $\alpha_2 < \alpha_1^\ast < x < y < u_1 < u_2$ and the following four relations become apparent:
\begin{align*}
u_1 y = y \alpha_2 x, \ \ \ \ \ \ \ \ \ \ \ \alpha_1^\ast y \alpha_2 = \alpha_2 u_2, \ \ \ \ \ \ \ \ \ \ \ \alpha_1^\ast u_1 = \alpha_2 x \alpha_1^\ast, \ \ \ \ \ \ \ \ \ \ \ u_2 x = x \alpha_1^\ast y.
\end{align*}
The only overlap $\alpha_1^\ast u_1y$ resolves as
\[
	\alpha_1^\ast (y \alpha_2 x) - (\alpha_2 x \alpha_1^\ast) y = \alpha_2 u_2 x - \alpha_2 x \alpha_1^\ast y = \alpha_2 x \alpha_1^\ast y - \alpha_2 x \alpha_1^\ast y = 0.
\]
A Gr\"obner basis calculation as in Example \ref{exp41} shows that this gives a presentation of the fixed ring $A^G$. For more details check the author's dissertation \cite[Example 4.2.3]{dissertation}.

In order to find the minimal projective resolutions we use the notation from before, i.e. $S_i = e_iA/(e_iA)_{\geq 1} = e_iA^G/(e_iA^G)_{\geq 1}$ and $P_i = e_iA^G$. The proof of \cite[Theorem 3.2]{MR2355031} again gives the beginnings of the easy minimal projective resolutions for $S_1$ and $S_3$:
\[
\begin{tikzcd}
0 \arrow[r] & P_2 \arrow[rr, "\begin{pmatrix} y \\ -\alpha_2 x \end{pmatrix} \boldsymbol{\cdot}"] & & P_1 \oplus P_2 \arrow[rr, "(u_1 \ \ y) \ \boldsymbol{\cdot}"] & & P_1\arrow[r] & S_1 \arrow[r] & 0, & \text{ and} \\
0 \arrow[r] & P_2 \arrow[rr, "\begin{pmatrix} \alpha_1^\ast y \\ -x \end{pmatrix} \boldsymbol{\cdot}"] & & P_2 \oplus P_3 \arrow[rr, "(x \ \ u_2) \ \boldsymbol{\cdot}"] & & P_3\arrow[r] & S_3 \arrow[r] & 0.
\end{tikzcd} 
\]
As in the previous example, exactness at the $0$-th and $1$-st position is by construction and left multiplication by $y$ ($x$, respectively) is injective due to the relations and the ordering. This shows we found the minimal projective resolutions of $S_1$ and $S_3$. Let's look at $S_2$:
\[
\begin{tikzcd}
\ldots \arrow[r] & P_1 \oplus P_3 \arrow[rrr, "\begin{pmatrix} \ \ u_1 \ \ \ \ y \alpha_2 \\ -x\alpha_1^\ast \ \ \ -u_2 \end{pmatrix} \boldsymbol{\cdot}"] & & & P_1 \oplus P_3 \arrow[rr, "(\alpha_1^\ast \ \ \alpha_2) \ \boldsymbol{\cdot}"] & & P_2 \arrow[r] & S_2 \arrow[r] & 0.
\end{tikzcd} 
\]
Again we have exactness in the $0$-th and $1$-st spot by \cite[Theorem 3.2]{MR2355031}. We need to understand the kernel of the map
\[
	\Gamma:	P_1 \oplus P_3 \longrightarrow P_1 \oplus P_3, \ \ \ \begin{pmatrix} \beta_1 \\ \beta_2 \end{pmatrix} \mapsto \begin{pmatrix} u_1 & y \alpha_2 \\ -x\alpha_1^\ast & -u_2	\end{pmatrix} \cdot \begin{pmatrix} \beta_1 \\ \beta_2 \end{pmatrix}.
\]
Let $\begin{pmatrix} f \\ g \end{pmatrix} \in \ker(\Gamma) \subseteq P_1 \oplus P_3$. Due to the relations, we can write $f = \sum\limits_{j=0}^{N}{c_ju_1^j} + c yf_1$ and $g = \sum\limits_{\ell = 0}^{M}{d_\ell u_2^\ell} + d x g_2$ for $c,c_j,d,d_\ell \in k$ as well as $M,N \in \NN_{0}$ and $f_1, \ g_2 \in P_2$. Then we obtain
\begin{align*}
&\begin{pmatrix} 0 \\ 0 \end{pmatrix} = \Gamma\left(\begin{pmatrix} f \\ g \end{pmatrix} \right) = \begin{pmatrix} u_1 & y \alpha_2 \\ -x\alpha_1^\ast & -u_2	\end{pmatrix} \cdot \begin{pmatrix} \sum\limits_{j=0}^{N}{c_ju_1^j} + c yf_1 \\ \sum\limits_{\ell = 0}^{M}{d_\ell u_2^\ell} + d x g_2 \end{pmatrix} \\
&\ \ \ = \begin{pmatrix} \sum\limits_{j=0}^{N}{c_ju_1^{j+1}} + c u_1yf_1 + \sum\limits_{\ell = 0}^{M}{d_\ell y\alpha_2u_2^\ell} + d y \alpha_2 x g_2 \\
-\sum\limits_{j=0}^{N}{c_jx \alpha_1^\ast u_1^j} - c x \alpha_1^\ast yf_1 - \sum\limits_{\ell = 0}^{M}{d_\ell u_2^{\ell+1}} - d u_2 x g_2 \end{pmatrix}.
\end{align*}
In the first entry, only $\sum\limits_{j=0}^{N}{c_ju_1^{j+1}}$ consists of only nonstar arrows while the rest of terms always contain $y$, which has two star arrows. This implies $c_j = 0$ for all $j$. Similarly, $\sum\limits_{\ell = 0}^{M}{d_\ell u_2^{\ell+1}}$ are the only summands with only star arrows in the second entry and so $d_\ell = 0$ for all $\ell$. Hence, we can assume that $f = y f_1$ and $g = x g_2$.

It is easy to check that $\begin{pmatrix} yf \\ -xf \end{pmatrix} \in \ker(\Gamma)$ for all $f \in P_2$:
\[
	\begin{pmatrix}	u_1 & y \alpha_2 \\ -x\alpha_1^\ast & -u_2 	\end{pmatrix} \cdot \begin{pmatrix} yf \\ -xf \end{pmatrix} = \begin{pmatrix} u_1yf - y\alpha_2xf \\ -x\alpha_1^\ast y f + u_2 x f \end{pmatrix}= \begin{pmatrix} (u_1y - y\alpha_2x)f \\ (-x\alpha_1^\ast y + u_2 x ) f \end{pmatrix} = 0.
\]
This implies that $\begin{pmatrix} 0 \\ x (f_1 + g_2) \end{pmatrix} \in \ker(\Gamma)$ as we know that $\begin{pmatrix} yf_1 \\ xg_2 \end{pmatrix}, \ \begin{pmatrix} y(-f_1) \\ xf_1 \end{pmatrix} \in \ker(\Gamma)$. Therefore, we can compute
\[
\begin{pmatrix}	u_1 & y \alpha_2 \\ -x\alpha_1^\ast & -u_2 	\end{pmatrix} \cdot \begin{pmatrix} 0 \\ x (f_1 + g_2) \end{pmatrix} = \begin{pmatrix} y \alpha_2 x (f_1 + g_2) \\ -u_2x (f_1 + g_2) \end{pmatrix} = \begin{pmatrix} 0 \\ 0 \end{pmatrix}.
\]
The first entry, $y\alpha_2 x (f_1 + g_2) = 0$, implies that $g_2 = -f_1$. This follows from $y \alpha_2 x$ not appearing as or including a leading term of one of the relations. To sum up, we can write
\[
	\ker(\Gamma) = \left\{\left. \begin{pmatrix} y \\ -x \end{pmatrix} \cdot f_1 \right| f_1 \in P_2\right\}
\]
and the minimal projective resolution becomes
\[
\begin{tikzcd}
0 \arrow[r] & P_2 \arrow[r, "\begin{pmatrix} y \\ -x \end{pmatrix} \boldsymbol{\cdot}"] & P_1 \oplus P_3 \arrow[rrr, "\begin{pmatrix} \ \ u_1 \ \ \ \ y \alpha_2 \\ -x\alpha_1^\ast \ \ \ -u_2 \end{pmatrix} \boldsymbol{\cdot}"] & & & P_1 \oplus P_3 \arrow[rr, "(\alpha_1^\ast \ \ \alpha_2) \ \boldsymbol{\cdot}"] & & P_2 \arrow[r] & S_2 \arrow[r] & 0,
\end{tikzcd} 
\]
which is injective at the last spot for the same reason as before: left multiplication by $x$ is injective due to the relations and the order. Once again, \cite[Proposition 3.19 (3)]{DanToDo} says the global dimension of $A^G$ equals the maximal projective dimension of the simple modules and thus $\gldim(A^G) = \pdim(S_2) = 3$.
\end{ex}
The last piece of this section helps understanding graded automorphisms $g$ of $\Pi_Q(k)$ for $Q = \widetilde{A_2}$ with $1$ as a pole of order $\GK(\Pi_Q(k)) = 2$. It was shown in \cite[Proposition 1.8]{MR2434290} that for Artin-Schelter regular domains, the denominator of the trace of a nonidentity automorphism cannot have $1$ as a root of order equal to the $\GKt$-dimension of the algebra. This result encouraged their definition of a quasi-reflection. Kirkman, Kuzmanovich and Zhang call a graded automorphism $\sigma$ of an AS-regular domain $S$ a quasi-reflection in \cite[Definition 2.2]{MR2434290} if it has $1$ as a root of the denominator of its trace $\Tr_S(\sigma,t) = p_\sigma(t)^{-1}$ of multiplicity equal to $\GK(S) - 1$. In contrast, for preprojective algebras many graded automorphisms have $1$ as a pole of the trace of order $2 = \GK(\Pi_Q(k))$. An easy example to see this occurs for $Q = \widetilde{A_2}$ with $c_1 = t_1 = 1$ and $c_2 = c_3 = t_2 = t_3 = -1$. Using the formula from the beginning of the section, we obtain $\Tr_A(g,t) = \dfrac{3-5t+3t^2}{(1-t)(1-t^3)}.$

A straight-forward approach to define quasi-reflections for preprojective algebras would be as follows. Instead of looking at $1$ as a root of the denominator of the trace function, one could consider the order of $1$ as a pole of $\Tr(g,t)$. The natural definition would be to call an automorphism $g$ a quasi-reflection if the order of $1$ as a pole of $\Tr(g,t)$ equals $\GKt$-dimension minus one. To get a feel for this idea one needs to understand the nontrivial graded automorphisms whose trace have $1$ as a pole of order $2$. As the following proposition shows in this case, at least for $n = 3$ the fixed rings $A^{\langle g \rangle}$ always have infinite global dimension which means we do not want to call these quasi-reflections.

\begin{prop}\label{prop44}
Let $g \neq \id$ be a graded automorphism of $A = \Pi_Q(k)$ for $Q = \widetilde{A_2}$. Denote $\Tr_A(g,t) = \frac{p(t)}{q(t)}$ and keep the formula given in the beginning of this section in mind. Then, the numerator $p(t)$ factors for $t = 1$ as follows:
\begin{align*}
p(1) &= (3 + c_1 + c_2 + c_3 + c_1 c_2 + c_2c_3 + c_1c_3)(1+(c_1t_1) + (c_1t_1)^2) \ \ \ \ \ \ \ \text{ if } c_1c_2c_3 = 1, \\
p(1) &= (3 + t_1 + t_2 + t_3 + t_1t_2 + t_2t_3 + t_1t_3)(1+(c_1t_1) + (c_1t_1)^2) \ \ \ \ \ \ \ \ \ \text{ if } t_1 t_2 t_3 = 1.
\end{align*}
If $\Tr_A(g,t)$ has $1$ as a pole of order $2$, then $A^{\langle g \rangle}$ has infinite global dimension.
\begin{pro}
Recall the formula of the trace due to Proposition \ref{prop34}:
\begin{small}
\begin{align*}
	\Tr_A(g,t) &= \frac{3 + (c_1+c_2+c_3+t_1+t_2+t_3)t + (c_1c_2 + c_1c_3 + c_2c_3 + t_1t_2 + t_1t_3 + t_2t_3 + 3c_1t_1)t^2}{(1-c_1c_2c_3t^3)(1-t_1t_2t_3t^3)} \\
	&\ \ \ + \frac{c_1t_1 (c_1 + c_2 + c_3 + t_1 + t_2 + t_3)t^3 + 3 (c_1t_1)^2t^4}{(1-c_1c_2c_3t^3)(1-t_1t_2t_3 t^3)}.
\end{align*}
\end{small}

Let $g$ be such that $\Tr_A(g,t)$ has $1$ as a pole of order $2$. The relations force $c_1 t_1 = c_2 t_2 = c_3t_3$. Also, in order to have $1$ as a pole of order $2$, we need $c_1c_2c_3 = t_1t_2t_3 = 1$. Hence,
\[
	(c_1t_1)^3 = c_1c_2c_3 \cdot t_1 t_2 t_3 = 1
\]
and so $c_1t_1$ must be a third root of unity. Consider the following calculation which shows how $p(1)$ can be factored
\begin{align*}
	&[3 + (c_1 + c_2 + c_3 + c_1c_2 + c_2c_3 + c_1c_3)] \cdot [1+c_1t_1 + (c_1t_1)^2] \\
	&\ \ \ = 3 + (c_1 + c_2 + c_3 + c_1c_2 + c_2c_3 + c_1c_3) + 3c_1t_1 + c_1t_1 (c_1 + c_2 + c_3 + c_1c_2 + c_2c_3 + c_1c_3) \nonumber \\
	& \ \ \ \ \ \ + 3(c_1t_1)^2 + (c_1t_1)^2(c_1 + c_2 + c_3 + c_1c_2 + c_2c_3 + c_1c_3) \nonumber \\[0.5em]
	&\ \ \ = 3 + (c_1 + c_2 + c_3 + (c_1t_1)(c_1c_2 + c_2c_3 + c_1c_3)) \nonumber \\
	&\ \ \ \ \ \ \ \ \ + (c_1c_2 + c_2c_3 + c_1c_3 + (c_1t_1)^2(c_1 + c_2 + c_3) + 3c_1t_1) \nonumber \\
	&\ \ \ \ \ \ \ \ \ + ((c_1t_1)(c_1 + c_2 + c_3) + (c_1t_1)^2(c_1c_2 + c_2c_3 + c_1 c_3)) + 3(c_1t_1)^2 \nonumber \\[0.5em]
	&\ \ \ = 3 + (c_1 + c_2 + c_3 + t_3 + t_1 + t_2) + (c_1c_2 + c_2c_3 + c_1c_3 + t_2t_3 + t_1t_3 + t_1t_2 + 3c_1t_1) \nonumber \\
	&\ \ \ \ \ \ \ \ \ + (c_1t_1)(c_1 + c_2 + c_3 + t_3 + t_1 + t_2) + 3(c_1t_1)^2. \nonumber
\end{align*}
This means if $1 + c_1t_1 + (c_1t_1)^2 = 0$ or $3 + (c_1 + c_2 + c_3 + c_1c_2 + c_2c_3 + c_1c_3) = 0$ we can have $1$ as a pole of order at most $1$. As $(c_1t_1)^3 = 1$, the only possibility to have $1 + c_1t_1 + (c_1t_1)^2$ not equal to zero is if $c_1t_1 = 1$. To sum up, our choices for the $c_i$ and $t_i$ are
\begin{center}
\begin{tabular}{lllllll}
	$\alpha_1 \mapsto c_1 \alpha_1$, & & & $\alpha_2 \mapsto c_2 \alpha_2$, & & & $\alpha_3 \mapsto \frac{1}{c_1c_2} \alpha_3$, \\[10pt]
	$\alpha_1^\ast \mapsto \frac{1}{c_1} \alpha_1^\ast$, & & & $\alpha_2^\ast \mapsto \frac{1}{c_2} \alpha_2^\ast$, & & & $\alpha_3^\ast \mapsto c_1c_2 \alpha_3^\ast$.	 
\end{tabular}
\end{center}
Since $c_1c_2c_3 = 1 = t_1t_2t_3$ the argument to find the generators of the fixed ring used in Example \ref{exp42} applies: it is enough to look at fixed paths of degree at most four as every fixed path can be decomposed as a product of fixed paths of degree less than or equal to four. We look at the following two cases.
\begin{itemize}[topsep=0pt]
	\item[(a)]	Either $c_1 = 1, \ c_2 \neq 1$, or $c_2 = 1, c_1 \neq 1$ or $c_1 \neq 1, c_2 \neq 1, c_1c_2 = 1$. The following list of elements helps to understand the fixed ring:
	\begin{center}
	\begin{tabular}{l|l}
	eigenvalue & eigenvector \\ \hline
	$1$ & $e_1, \ e_2, \ e_3, \ \alpha_1 \alpha_1^\ast, \ \alpha_2 \alpha_2^\ast, \ \alpha_3 \alpha_3^\ast, \ \alpha_1 \alpha_2 \alpha_3, \ \alpha_3^\ast \alpha_2^\ast \alpha_1^\ast, \ \alpha_2 \alpha_3 \alpha_1, \ \alpha_1^\ast \alpha_3^\ast \alpha_2^\ast,$ \\
	& $\alpha_3 \alpha_1 \alpha_2, \ \alpha_2^\ast \alpha_1^\ast \alpha_3^\ast, \ \alpha_1 \alpha_2 \alpha_2^\ast \alpha_1^\ast, \ \alpha_2 \alpha_3 \alpha_3^\ast \alpha_2^\ast, \ \alpha_3 \alpha_1 \alpha_1^\ast \alpha_3^\ast$ \\ \hline
	$c_1$ & $\alpha_1, \ \alpha_3^\ast \alpha_2^\ast, \ \alpha_1 \alpha_2 \alpha_2^\ast, \ \alpha_1 \alpha_2 \alpha_3 \alpha_1, \ \alpha_1 \alpha_1^\ast \alpha_3^\ast \alpha_2^\ast $ \\ \hline 
	$c_1^{-1}$ & $\alpha_1^\ast, \ \alpha_2 \alpha_3, \ \alpha_2 \alpha_2^\ast \alpha_1^\ast, \ \alpha_2 \alpha_3 \alpha_1 \alpha_1^\ast, \ \alpha_1^\ast \alpha_3^\ast \alpha_2^\ast \alpha_1^\ast$ \\ \hline
	$c_2$ & $\alpha_2, \ \alpha_1^\ast \alpha_3^\ast, \ \alpha_2\alpha_3 \alpha_3^\ast, \ \alpha_2\alpha_3 \alpha_1 \alpha_2, \ \alpha_2 \alpha_2^\ast \alpha_1^\ast \alpha_3^\ast$ \\ \hline
	
	$c_2^{-1}$ & $\alpha_2^\ast, \ \alpha_3 \alpha_1, \ \alpha_3 \alpha_3^\ast \alpha_2^\ast, \ \alpha_3 \alpha_1 \alpha_2 \alpha_2^\ast, \ \alpha_2^\ast \alpha_1^\ast \alpha_3^\ast \alpha_2^\ast$ \\ \hline
	
	$c_1c_2$ & $\alpha_3^\ast, \ \alpha_1\alpha_2, \ \alpha_1\alpha_1^\ast \alpha_3^\ast, \ \alpha_1\alpha_2 \alpha_3 \alpha_3^\ast, \ \alpha_3^\ast \alpha_2^\ast \alpha_1^\ast \alpha_3^\ast$ \\ \hline
	
	$(c_1c_2)^{-1}$ & $\alpha_3, \ \alpha_2^\ast \alpha_1^\ast, \ \alpha_3\alpha_1 \alpha_1^\ast,\ \alpha_3 \alpha_1 \alpha_2 \alpha_3, \ \alpha_3 \alpha_3^\ast \alpha_2^\ast \alpha_1^\ast$
	\end{tabular}
	\end{center}
	So, if precisely one of $c_1, \ c_2, \ c_1c_2$ equals $1$, we obtain the same diagram up to a rotation. For example in case $c_1c_2 = 1$, we get the path algebra underlying $A^G$ as
	
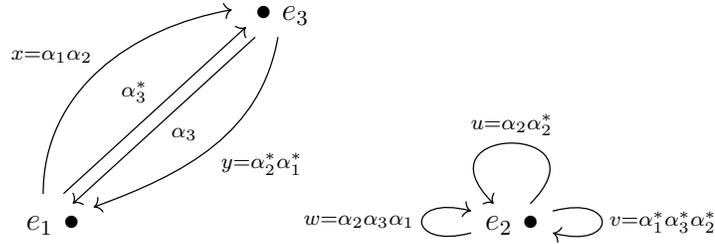
\[
\begin{tikzcd}
& & & & \bullet \ e_3 \arrow[llddd,"\alpha_3"] \arrow[llddd,"y =\alpha_2^\ast \alpha_1^\ast",out=260,in=20,shift left = 0ex] \\
& \\
& \\
& & e_1 \ \bullet \arrow[rruuu,"\alpha_3^\ast", shift left = 1ex] \arrow[rruuu,"x = \alpha_1 \alpha_2",out = 90, in = 190,shift left = 1ex] & & & & e_2 \ \bullet \arrow[loop right, "v = \alpha_1^\ast \alpha_3^\ast \alpha_2^\ast"] \arrow[loop left, "w =\alpha_2 \alpha_3 \alpha_1"] \arrow[swap,loop, "u = \alpha_2 \alpha_2^\ast"]
\end{tikzcd}
	\] 
	with relations $uw = vu, \ uv = vu, \ vw = wv$ and $wv = u^3$ for the right connected component and relations $xy = (\alpha_3^\ast \alpha_3)^2, \ yx = (\alpha_3 \alpha_3^\ast)^2, \ x \alpha_3 \alpha_3^\ast = \alpha_3^\ast \alpha_3 x, \ y\alpha_3^\ast \alpha_3 = \alpha_3 \alpha_3^\ast y$ for the left connected component. Note that $\alpha_1 \alpha_1^\ast = \alpha_3^\ast \alpha_3$ and $\alpha_2 \alpha_3 \alpha_3^\ast \alpha_2^\ast = u^2$. \\
	
Consider the right component $k[u,v,w]/(u^3 - vw)$ as the coordinate ring of the affine variety $Y = Z(u^3 - vw)$. The beginning of \cite[Chapter 5]{MR0463157} says that the point $P = (0,0,0)$ is singular since the Jacobian matrix $\begin{pmatrix} 3u^2 & -w & -v \end{pmatrix}$ is zero at $P$. Then, \cite[Theorem 5.1]{MR0463157} implies that the local ring $\mathcal{O}_{P,Y}$ is not regular and as a localization this also makes the local ring $k[u,v,w]/(u^3 - vw)$ not regular. By \cite[Th\'eor\`eme 3]{MR0086071}, a noetherian local ring which is not regular must have infinite global dimension.
	
	\item[(b)]	All three of $c_1, \ c_2, \ c_1c_2$ are not equal to $1$. Hence, the only elements fixed of degree at most four are 
	\begin{align*}
		&e_1, \ e_2, \ e_3, \\
		&\alpha_1 \alpha_1^\ast, \ \alpha_2 \alpha_2^\ast, \ \alpha_3 \alpha_3^\ast, \\
		&\alpha_1 \alpha_2 \alpha_3, \ \alpha_2 \alpha_3 \alpha_1, \ \alpha_3 \alpha_1 \alpha_2, \ \alpha_3^\ast \alpha_2^\ast \alpha_1^\ast, \ \alpha_1^\ast \alpha_3^\ast \alpha_2^\ast, \ \alpha_2^\ast \alpha_1^\ast \alpha_3^\ast, \\
		&\alpha_1 \alpha_2 \alpha_2^\ast \alpha_1^\ast, \ \alpha_2 \alpha_3 \alpha_3^\ast \alpha_2^\ast, \ \alpha_3 \alpha_1 \alpha_1^\ast \alpha_3^\ast.
	\end{align*}
Therefore, the fixed ring is isomorphic to the path algebra
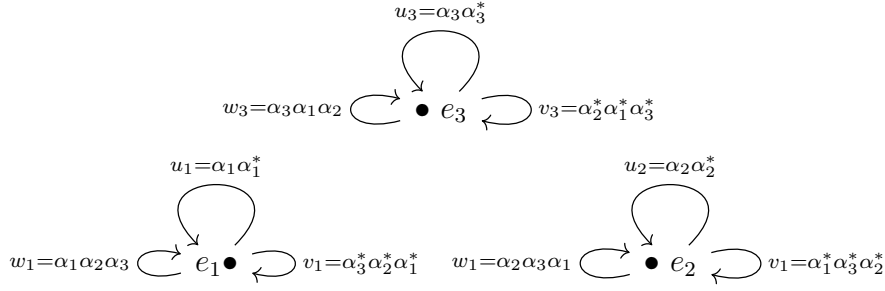
\[
\begin{tikzcd}
& & & & \bullet \ e_{3} \arrow[swap,loop, "u_3 = \alpha_3 \alpha_3^\ast"] \arrow[loop right, "v_3 = \alpha_2^\ast \alpha_1^\ast \alpha_3^\ast"] \arrow[loop left, "w_3 =\alpha_3 \alpha_1 \alpha_2"] \\
& \\
& & e_1 \bullet \arrow[swap,loop, "u_1 = \alpha_1 \alpha_1^\ast"] \arrow[loop right, "v_1 = \alpha_3^\ast \alpha_2^\ast \alpha_1^\ast"] \arrow[loop left, "w_1 =\alpha_1 \alpha_2 \alpha_3"] & & & & \bullet \ e_2 \arrow[swap,loop, "u_2 = \alpha_2 \alpha_2^\ast"] \arrow[loop right, "v_1 = \alpha_1^\ast \alpha_3^\ast \alpha_2^\ast"] \arrow[loop left, "w_1 =\alpha_2 \alpha_3 \alpha_1"]
\end{tikzcd}
\] 
with relations $u_iw_i = w_iu_i, \ u_iv_i = v_iu_i, \ v_iw_i = w_iv_i$ and $w_iv_i = u_i^3$ for $i = 1,2,3$. The same reasoning as in part (a) proves that this algebra has infinite global dimension.\qed
\end{itemize}
\end{pro}
\end{prop}

\section{Properties of the Fixed Ring of a non-connected Algebra}\label{sectiongeneral}
In this section, let $B$ be a noetherian $\NN$-graded twisted Calabi-Yau algebra of global dimension $d \geq 1$ with degree zero piece $B_0 \cong k^n$. Lemma \ref{lemma27} combined with Proposition \ref{prop29} implies that $B$ and $B^{\op}$ satisfy the $\chi$-condition. Denote the primitive orthogonal idempotents by $e_i$. This gives us all simple right $B$-modules $S_i = e_iB/(e_iB)_{\geq 1}$ and all simple left $B$-modules $S_i^\ast = Be_i / (Be_i)_{\geq 1}$, respectively. Fix a finite subgroup $G$ of the graded automorphism group of $B$ such that every $g \in G$ satisfies $g(e_i) = e_i$. Notice that this implies $(B^G)_0 = B_0 \cong k^n$.

We will present connections between $B^G$ having finite global dimension together with satisfying the generalized Gorenstein condition and $\gldim(B) = \gldim(B^G)$ as well as $B$ being a projective $B^G$-module. In particular, the only direction missing from proving that the mentioned three conditions are equivalent is to show $\gldim(B^G) = \gldim(B)$ implies $B^G$ satisfies the generalized Gorenstein condition. We will be able to obtain this remaining piece under some stronger assumptions in the next section.

In particular, if $Q$ is the quiver corresponding to an extended Dynkin graph of type $\widetilde{A_{n-1}}, \ \widetilde{D_{n-1}}$ or $\widetilde{E_m}$ for $m = 6,7,8$ then $A = \Pi_Q(k)$ is an $\NN$-graded Calabi-Yau algebra of global dimension $2$ whose degree zero piece is isomorphic to $k^n$ (where $n = m+1$ if $Q = \widetilde{E_m}$). Moreover, $A$ is noetherian and satisfies the $\chi$-condition by Proposition \ref{prop210}. Thus, for $G$ a nontrivial finite subgroup of graded automorphisms of $A$ satisfying that every $g \in G$ fixes the primitive idempotents $e_i$ for $i = 1,\ldots,n$ all results of this section apply.

Before we can start we need to introduce more notation. Most of the conditions are required in order to apply theorems from \cite{MR1304753}. Therefore, let's adopt the notation from \cite{MR1304753}: $\qgr(B)$ denotes the \textit{quotient category} $\gr(B)/\tors(B)$ where $\gr(B)$ is the category of all finitely generated graded right $B$-modules and $\tors(B)$ denotes the subcategory of finite dimensional modules. For every module $M \in \gr(B)$, we denote its image in $\qgr(B)$ by $\mathcal{M} = \pi(M)$.

Also recall the definition of \textit{cohomological dimension of $\proj(B)$} from \cite[p.276]{MR1304753} which uses $H^i(\mathcal{M}) = H^i_B(\mathcal{M}) := \Ext^i_{\qgr(B)}(\mathcal{B},\mathcal{M})$:
\[
\cd(\proj B) = \max\{i \mid H^i(\mathcal{M}) \neq 0 \text{ for some } \mathcal{M} \in \qgr(B)\}.
\]

Since $\gldim(B) = d$ is finite, \cite[Proposition 7.10(3)]{MR1304753} says $\cd(\proj B) \leq d-1$ is finite as well. Then, we can use \cite[Proposition 7.10(1)]{MR1304753} to get a different description of $\cd(\proj B)$ involving $\underline{H}^i(\mathcal{M}) := \bigoplus\limits_{s = - \infty}^{\infty}{H^i(\mathcal{M}[s])}$:
\[
	\cd(\proj B) = \max\{i \mid \underline{H}^i(\mathcal{B}) \neq 0\}.
\]
This says, we only need to look at $\mathcal{B}$ in $\qgr(B)$. Even better, for $i \geq 1$, we have the equality $\underline{H}^i(\mathcal{B}) = \lim\limits_{m \to \infty} \Ext^{i+1}_B(B/B_{\geq m},B)$ thanks to \cite[Proposition 7.2(2)]{MR1304753}. Moreover, the same proposition says $\underline{H}^0(\mathcal{B})$ surjects onto $\lim\limits_{m \to \infty} \Ext^1_B(B/B_{\geq m},B)$.

The following lemma connects the cohomological dimension to the global dimension of $B$.

\begin{lemma}\label{lemma51}
$\cd(\proj B) = \gldim(B)-1$.
\begin{pro}
We obtain $\cd(\proj B) \leq \gldim(B) - 1$ from \cite[Proposition 7.10(3)]{MR1304753}. To achieve the equality $\cd(\proj B) = \gldim(B) - 1$, we consider the following constructions and results to generalize the idea of \cite[Theorem 8.1(3)]{MR1304753}. For $d = \gldim(B) \geq 1$ we are going to show that $\lim\limits_{m \to \infty}\Ext^{d}_B(B/B_{\geq m},B) \neq 0$. By our assumptions $B$ is twisted Calabi-Yau. This guarantees the existence of the $k$-central invertible bimodule $U = \Ext_{B^e}^d(B,B^e)$ as in Definition \ref{CYdef}. From \cite[Proposition 4.12]{DanToDo} it follows that
\begin{align*}
\lim_{m \to \infty}{\Ext^d_B(B/B_{\geq m},B)} &\cong \lim_{m \to \infty} \Tor_{d-d}^B\left((B/B_{\geq m})^\ast \otimes_B U , \ B\right) \\
&\cong \lim_{m \to \infty} \left((B/B_{\geq m})^\ast \otimes_B U \right) \otimes_B B \\
&\cong \lim_{m \to \infty} (B/B_{\geq m})^\ast \otimes_B U. \\
\end{align*}
We can calculate the direct limit of $(B/B_{\geq m})^\ast$. Applying the left-exact contravariant functor $\Hom(-,k)$ to the natural projections $B/B_{\geq m+1} \longrightarrow B/B_{\geq m} \longrightarrow 0$, we obtain exact maps
\[
\begin{tikzcd}
0 \arrow[r] & \Hom_k(B/B_{\geq m},k) \arrow[r] & \Hom_k(B/B_{\geq m+1},k).
\end{tikzcd}
\]
Hence, the direct limit becomes the union of all $\Hom_k(B/B_{\geq m},k)$ which is nonzero. This does not change if we apply $\otimes_B U$. 

Thus, for $d = \gldim(B) \geq 2$, we showed that $\underline{H}^{d-1}(\mathcal{B}) = \lim\limits_{m \to \infty}{\Ext^d_B(B/B_{\geq m},B)} \neq 0$ using \cite[Proposition 7.2(2)]{MR1304753}. For $d = \gldim(B) = 1$, we have the surjection of $\underline{H}^0(\mathcal{B})$ onto $\lim\limits_{m \to \infty} \Ext^1_B(B/B_{\geq m},B) \neq 0$. Together, for $\gldim(B) \geq 1$ we obtain $\underline{H}^{d-1}(\mathcal{B}) \neq 0$. Hence $\cd(\proj B) = d-1$.\qed
\end{pro}
\end{lemma}

\begin{prop}\label{prop52}
If $B^G$ is twisted Calabi-Yau, then $\gldim(B^G) = \gldim(B)$.
\begin{pro}
The proof of \cite[Lemma 1.11(d)]{MR2434290} works in our case. We immediately gain that $B^G$ is a noetherian $\NN$-graded algebra (see Lemma \ref{lemma214}(b)). Moreover, $B^G$ satisfies the $\chi$-condition due to Proposition \ref{prop217}. These properties will be important to use results from \cite{MR1304753}. 

Now, Lemma \ref{lemma51} applies to both $B$ and $B^G$ saying $\cd(\proj B) = \gldim(B) -1$ as well as $\cd(\proj B^G) = \gldim(B^G) - 1$. Therefore, in order to finish the proof, it suffices to show that $\cd(\proj B) = \cd(\proj B^G)$.

Recall that $B$ satisfies the $\chi$-condition and the natural inclusion of $B^G \hookrightarrow B$ makes $B$ a finitely generated $B^G$-module by Proposition \ref{prop216}. Hence, \cite[Corollary 8.4(1)]{MR1304753} says $\cd(\proj B) \leq \cd(\proj B^G)$. Further, the natural projection
\[
	\pi: B \to B^G, \ b \mapsto \frac{1}{|G|}\sum_{g \in G}{g(b)}
\]
allows one to decompose $B = B^G \oplus C$ as $(B^G, B^G)$-bimodules (see Lemma \ref{lemma214}(a)). 

The same calculation as in \cite[Lemma 1.11(d)]{MR2434290} for $i = d-1$ provides
\begin{align*}
	\cd(\proj B) &= \max\{i \mid \underline{H}^i_B(\mathcal{B}) \neq 0\} \\
&= \max\{i \mid \underline{H}^i_B(\mathcal{B} \oplus \mathcal{C}) \neq 0\} \\	
&= \max\{i \mid \underline{H}^i_{B^G}(\mathcal{B} \oplus \mathcal{C}) \neq 0\} \hspace{1cm} \text{\cite[Theorem 8.3(3)]{MR1304753}} \\	
&\geq \max\{i \mid \underline{H}^i_B(\mathcal{B^G}) \neq 0\} \\	
&= \cd(\proj B^G).
\end{align*}
This shows that $\gldim(B^G) - 1 = \cd(\proj B^G) = \cd(\proj B) = \gldim(B) -1$.
\qed
\end{pro}
\end{prop}

In the case of a noetherian Artin-Schelter regular algebra $S$ it cannot happen that $S^G$ has finite global dimension but $S^G$ is not AS-Gorenstein: Lemma \ref{lemma214} and Proposition \ref{prop216} apply and so whenever the fixed ring has finite global dimension, the Gorenstein condition comes along with it due to \cite[Lemma 1.10(c)]{MR2434290}. However, for preprojective algebras Example \ref{exp43} shows that $A^G$ having finite global dimension does not guarantee the generalized Gorenstein condition. 

Next, an easy lemma gives that $\gldim(B^G)$ is always greater than or equal to $\gldim(B)$. Afterwards, we prove that we obtain equality provided $B$ is a projective $B^G$-module.

\begin{lemma}\label{lemma53}
The fixed ring $B^G$ satisfies $\gldim(B^G) \geq \gldim(B)$.
\begin{pro}
If $\gldim(B^G) = \infty$ there is nothing to prove. So assume $\gldim(B^G)$ is finite. As in the proof of Proposition \ref{prop52}, $B$ and $B^G$ are right noetherian $\NN$-graded algebras. Since both $B$ and $B^G$ have finite global dimension, \cite[Proposition 7.10(3)]{MR1304753} gives $\cd(\proj B) \leq \gldim(B) - 1$ and $\cd(\proj B^G) \leq \gldim(B^G) - 1$. Since $B$ is twisted Calabi-Yau, Lemma \ref{lemma51} gives $\cd(\proj B) = \gldim(B) - 1$. Since $B$ satisfies the $\chi$-condition, \cite[Corollary 8.4(1)]{MR1304753} applied to $B^G \hookrightarrow B$ says $\cd(\proj B) \leq \cd(\proj B^G)$. Combining these results yields
\[
	\gldim(B) - 1 = \cd(\proj B) \leq \cd(\proj B^G) \leq \gldim(B^G) - 1
\]
which gives $\gldim(B) \leq \gldim(B^G)$ as claimed.\qed
\end{pro}
\end{lemma}

\begin{prop}\label{prop54}
If $B$ is a projective $B^G$-module then $\gldim(B) = \gldim(B^G)$. 
\begin{pro}
Lemma \ref{lemma53} provides $d = \gldim(B) \leq \gldim(B^G)$. The other inequality follows from the fact that a minimal projective resolution of $S_j = e_jB/e_jB_{\geq 1} = e_jB^G/(e_jB^G)_{\geq 1}$ as a $B$-module yields a projective resolution as a $B^G$-module. Hence, \cite[Proposition 3.19 (3)]{DanToDo} which says $\gldim(B^G) = \max\{\pdim(S_j) \mid j = 1,\ldots,n\}$ finishes the proof.\qed
\end{pro}
\end{prop}

After giving two necessary conditions for $\gldim(B) = \gldim(B^G)$, as a conclusion, we finish by showing that one of the conditions is actually equivalent.
\begin{prop}\label{prop55}
$B$ is a projective $B^G$-module if and only if $\gldim(B) = \gldim(B^G)$.
\begin{pro}
One direction is Proposition \ref{prop54}. Assume that $\gldim(B) = \gldim(B^G)$. By \cite[Proposition 3.19 (1)]{DanToDo}, we get $\pdim_{B^G}(B) = \sup\{q \in \NN \mid \Tor_q^{B^G}(B_0^G,B) \neq 0\}$ and so we need to prove that $\Tor_q^{B^G}(B_0^G,B) = 0$ for all $q \geq 1$. By \cite[Theorem 5.6.6]{MR1269324} applied to $B^G \hookrightarrow B$ with $M = B_0^G$ and some left $B^G$-module $N$ we obtain the spectral sequence
\begin{align}\label{eqweibel}
	E_{pq}^2 = \Tor_p^B\left(\Tor_q^{B^G}(B_0^G,B),N\right) \R \Tor_{p+q}^{B^G}(B_0^G,N).
\end{align}
Thanks to the hypothesis, we can denote $d = \gldim(B) = \gldim(B^G)$ as well as call $e = \pdim_{B^G}(B) = \max \{q \in \NN \mid {\Tor_q^{B^G}(B_0^G,B) \neq 0}\} \leq d$. We need to show that $e = 0$. By using a minimal projective resolution of $B$ as $B^G$-module of length $e = \pdim_{B^G}(B)$, it follows that $\Tor_q^{B^G}(B_0^G,B) = 0$ for all $q > e$. Analogously, $\Tor_p^B(-,-) = 0$ for $p > d$. Therefore, the only nonzero terms of $E_{pq}^2$ appear for $0 \leq p \leq d$ and $0 \leq q \leq e$ making $E_{de}^2$ the top right corner of this rectangle. Since we know that $E_{pq}^2 \R \Tor_{p+q}^{B^G}(B_0^G,N)$ we need to understand $E_{pq}^r$ for $r \geq 2$. The data includes the differentials $d_{pq}^r: E_{pq}^r \to E_{p-r,q+r-1}$ and thus for $r \geq 2$ taking homology at $E_{de}^r$ does not change $E_{de}^r$ due to the involved maps being trivial. Therefore, the convergence in \eqref{eqweibel} yields
\begin{align}\label{eq7}
	\Tor_d^B(\underbrace{\Tor_e^{B^G}(B_0^G,{}_{B^G}B_B)}_{= X_B},N) \cong \Tor_{d+e}^{B^G}(B_0^G,N).
\end{align}
First, we notice that $X$ is finite dimensional over $k$ since applying $B_0^G \otimes -$ to a minimal finitely generated projective resolution $P_\bullet \to 0$ of $B$ only leaves finite dimensional pieces at each step. This does not change if $X_B$ is calculated using a resolution of $B_0^G$ instead. Next, we prove that $\pdim(X) = \gldim(B)$. Being a $B$-module gives $\pdim(X) \leq \gldim(B)$ for free. For the converse, notice that every finite dimensional module over $B$ has a composition series with factors which are shifted copies of simple modules. Therefore, by short exact sequences and the generalized Gorenstein condition which $B$ satisfies by Proposition \ref{prop29}, we gain $\Ext_B^d(X,B) \neq 0$. This is enough to guarantee $\pdim(X) = d.$

From the characterization of $\pdim(X)$ and symmetry in the argument of \cite[Proposition 3.19 (1)]{DanToDo} it follows that $d = \gldim(B) = \pdim(X) = \max\{q \in \NN\mid \Tor_q^B(X,B_0) \neq 0\}$ where $B_0 = B/J(B)$. Therefore by choosing $N = B_0$, the left hand side of Equation \eqref{eq7} is nonzero. The only way the right hand side can be nonzero is if $d+e = d$ or $e = 0$. \qed
\end{pro}
\end{prop}
Proposition \ref{prop55} shows that $\gldim(B^G) = \gldim(B)$ if and only if $B$ is projective as a $B^G$-module. This is true for the preprojective algebra $A$. In particular, it finally finishes the proof of Example \ref{exp42} and makes sure that $A$ is a projective $A^G$-module. However, it can happen that $\gldim(A^G)$ is finite but different from $\gldim(A)$ as described in Example \ref{exp43}. We have now proved that in Example \ref{exp43} $A$ is not projective as an $A^G$-module. The situation for an Artin-Schelter regular algebra $S$ with a finite group $G$ acting on $S$ is nicer. In \cite[Lemma 1.10]{MR2434290} and \cite[Lemma 1.11]{MR2434290} it was proved that the fixed ring $S^G$ has finite global dimension if and only if $\gldim(S) = \gldim(S^G)$ if and only if the fixed ring is Artin-Schelter regular if and only if $S$ is free over $S^G$. 

\section{The main theorem}\label{section6}
Let $B$ be a noetherian $\NN$-graded algebra of finite global dimension $d \geq 1$ with degree zero piece $B_0 \cong k^n$. Denote the primitive orthogonal idempotents by $\{e_1,\ldots,e_n\}$. This gives us all simple right $B$-modules $S_i = e_iB/(e_iB)_{\geq 1}$ and all simple left $B$-modules $S_i^\ast = Be_i / (Be_i)_{\geq 1}$, respectively. Fix a finite subgroup $G$ of the graded automorphism group of $B$ such that every $g \in G$ satisfies $g(e_i) = e_i$. Notice that this implies $(B^G)_0 = B_0 \cong k^n$.

The purpose of this section is to find a converse of Proposition \ref{prop52}. In light of Example \ref{exp43} it became apparent that in order to achieve the generalized Gorenstein condition, all simple modules need to have the same projective dimension over $B^G$. More precisely, this section starts with technical results leading to Theorem \ref{thm63} which gives sufficient conditions for $B$ to satisfy the generalized Gorenstein condition. This is enough to prove the main result: if additionally $B$ is twisted Calabi-Yau and $\pdim_{B^G}(S_i) = \pdim_{B^G}(S_j)$ as well as $\pdim_{B^G}(S_i^\ast) = \pdim_{B^G}(S_j^\ast)$ for the simple modules $S_i = e_iB^G/(e_iB^G)_{\geq 1} = e_iB/(e_iB)_{\geq 1}$ and $S_i^\ast = B^Ge_i/(B^Ge_i)_{\geq 1} = Be_i/(Be_i)_{\geq 1}$ and all $i,j$, then the following are equivalent:
\begin{itemize}[topsep=0pt] \vspace{-0.5em}
	\item[(1)]	$\gldim(B^G) = \gldim(B)$.
	\item[(2)]	$B^G$ is twisted Calabi-Yau.
	\item[(3)]	$B$ is a projective $B^G$-module.
\end{itemize}

But first, let's see how the $\Ext$-groups involving the simple modules determine $\Ext$-groups for finite dimensional modules. The proof of this standard result using the long exact $\Ext$ sequence is omitted.

\begin{lemma}\label{lemma61}
If $F_B$ is a finite dimensional $B$-module and $\Ext^i_B(S_j,B) = 0$ for all $j$ and all $i \leq p$, then $\Ext^i_B(F,B) = 0$ for $i \leq p$.
\end{lemma}

\begin{prop}\label{prop62}
Let $B^{\op}$ be the opposite ring of $B$. Denote the global dimension of $B$ by $d = \gldim(B) < \infty$. Assume $\pdim_B(S_i) = \pdim_B(S_j)$ and $\pdim_{B^{\op}}(S_i^\ast) = \pdim_{B^{\op}}(S_j^\ast)$ for all $i,j$ where $S_i$ and $S_i^\ast$ denote the simple right and left $B$-modules. If \vspace{-.5em}
\begin{itemize}[topsep=0pt]
	\item[$\bullet$] $\Ext_{B^{\op}}^{i}(S_j^\ast,B) = 0 = \Ext_B^i(S_j,B)$ for $i \neq d$ and
	\item[$\bullet$] $\dim_k(\Ext_{B^{\op}}^d(S_j^\ast,B)) < \infty$ and $\dim_k(\Ext_B^d(S_j,B)) < \infty$ for all $j$
\end{itemize} \vspace{-0.5em} then $\dim_k(\Ext_B^d(S_j,B)) = 1$.
\begin{pro} 
This uses ideas from \cite[Proposition 3.1]{MR1371143}. Since $B$ is noetherian, \cite[Theorem 8.27]{MR2455920} implies $\gldim(B) = \gldim(B^{\op})$ and so $\pdim(S_j) = \pdim(S_1) = \pdim(S_1^\ast) = \pdim(S_j^\ast) = d$ for all $j$ by \cite[Proposition 3.19(3)]{DanToDo}. Consider a minimal projective resolution of $S_j$:
\[
\begin{tikzcd}
0 \arrow[r] & Q_d^{(j)} \arrow[r] & \ldots \arrow[r] & Q_2^{(j)} \arrow[r] & Q_1^{(j)} \arrow[r] & e_jB \arrow[r] & S_j \arrow[r] & 0,
\end{tikzcd}
\]
where each $Q_i^{(j)}$ equals a sum of shifted copies of the indecomposable projectives $e_mB[-\ell]$. If we apply $\Hom(-,B)$, Lemma \ref{lemma25} applies and the result is
\[
\begin{tikzcd}
0 & \arrow[l] \left(Q_d^{(j)}\right)^\ast & \arrow[l] \ldots & \arrow[l] \left(Q_2^{(j)}\right)^\ast & \arrow[l] \left(Q_1^{(j)}\right)^\ast & \arrow[l] Be_j & \arrow[l] 0,
\end{tikzcd}
\]
where $\left(Q_i^{(j)}\right)^\ast$ is the left $B$-module consisting of a sum of the corresponding indecomposable projectives $Be_m[\ell]$. Since $\Ext_B^i(S_j,B) = 0$ for $i \neq d$, the second resolution is exact everywhere but at the $d$-th position. Hence, it can be understood as a minimal projective resolution of $F_j^\ast = \Ext_B^d(S_j,B)$. Notice, the $\ast$ indicates that it is a left $B$-module. Since $\Hom(\Hom(e_mB,B),B) = e_mB$, it follows that $\Ext^d_{B^{\op}}(F_j^\ast,B) = S_j$ and $\Ext^i_{B^{\op}}(F_j^\ast,B) = 0$ for $i \neq d$. In particular, $\dim_k(\Ext^d_{B^{\op}}(F_j^\ast,B)) = 1$.

Define $w_j = \dim_k(\Ext^d_{B^{\op}}(S_j^\ast,B)) < \infty$. Notice, that $\gldim(B^{\op}) = \pdim(S_j^\ast) = d$. Therefore, we can take a minimal projective resolution $\mathcal{P}_{\bullet}$ of $S_j^\ast$ whose boundary maps can be described by matrices with entries of degree bigger than or equal to $1$ due to being minimal. After applying $\Hom_{B^{\op}}(-,B)$ to $\mathcal{P}_{\bullet}$ we can use a symmetric version of Lemma \ref{lemma25}. Now, at the $d$-th position we have
\[
\begin{tikzcd}
0 & \arrow[l] \Hom_{B^{\op}}(\mathcal{P}_d,B) & \arrow[l, swap, "\delta_d^\ast"] \ldots 
\end{tikzcd}
\]
where $\delta_d^\ast$ is a matrix with entries in $B^{\op}_{\geq 1}$. Clearly, taking homology at the described spot leaves behind the lowest degree part of $ \Hom_{B^{\op}}(\mathcal{P}_d,B)$ which implies $w_j > 0$. Let $M^\ast$ be a finite dimensional left $B$-module. We are going to show by induction on $\dim_k(M^\ast)$ that 
\[
	\dim_k(\Ext^d_{B^{\op}}(M^\ast,B)) = \sum\limits_{s = 1}^{\dim_k(M^\ast)}{w_{m_s}}, \ \ \ \text{ for some numbers } m_1, \ldots m_{\dim_k(M^\ast)} \in \{1,\ldots,n\}.
\]
If $\dim_k(M^\ast) = 1$ then $\Ext^d_{B^{\op}}(M^\ast,B) = \Ext^d_{B^{\op}}(S_{m_1}^\ast,B)$ and $\dim_k(\Ext^d_{B^{\op}}(S_{m_1}^\ast,B)) = w_{m_1}$. Otherwise, consider the short exact sequence
\[
\begin{tikzcd}
0 \arrow[r] & K^\ast \arrow[r] & M^\ast \arrow[r] & S_{m}^\ast \arrow[r] & 0 
\end{tikzcd}
\]
and notice that $\dim_k(K^\ast) = \dim_k(M^\ast) - 1$. Calculating the long exact $\Ext$ sequence gives
\[
\begin{tikzcd}
\Ext^{d-1}_{B^{\op}}(K^\ast,B) \arrow[r] & \Ext^d_{B^{\op}}(S_m^\ast,B) \arrow[r] & \Ext^d_{B^{\op}}(M^\ast,B) \arrow[r] & \Ext^d_{B^{\op}}(K^\ast,B) \arrow[r] & 0.
\end{tikzcd}
\]
Due to $\Ext^{r}_{B^{\op}}(S_j^\ast,B) = 0$ for all $r \leq d-1$, a symmetric version of Lemma \ref{lemma61} says $\Ext^{d-1}_{B^{\op}}(K^\ast,B)$ vanishes. Hence, the last three terms form a short exact sequence and the induction hypothesis applies to give
\begin{align*}
	\dim_k(\Ext^d_{B^{\op}}(M^\ast,B)) &= \dim_k(\Ext^d_{B^{\op}}(K^\ast,B)) + \dim_k(\Ext^d_{B^{\op}}(S_m^\ast,B)) \\
	&= \sum_{s = 1}^{\dim_k(M^\ast) - 1}{w_{m_s}} + w_{m_{\dim_k(M^\ast)}} = \sum_{s=1}^{\dim_k(M^\ast)}{w_{m_s}}.
\end{align*}
To sum up, we apply this result to $M^\ast = F_j^\ast = \Ext_B^d(S_j,B)$ which is finite dimensional by our assumption. This yields 
\[
	1 = \dim_k(\Ext^d_{B^{\op}}(F_j^\ast,B)) = \sum\limits_{s = 1}^{\dim_k(F_j^\ast)}{w_{m_s}}.
\]
Here, a sum of positive integers equals $1$ which can only happen with just one summand. Therefore, $\dim_k(\Ext_B^d(S_j,B)) = \dim_k(F_j^\ast) = 1$.\qed
\end{pro}
\end{prop}

The following theorem is a modification of \cite[Theorem 3.6, (2) $\R$ (1)]{MR2526379}.
\begin{thm}\label{thm63}
Let $B^{\op}$ be the opposite ring of $B$. Suppose that both $B$ and $B^{\op}$ satisfy the $\chi$-condition. Assume that $\pdim_B(S_i) = \pdim_B(S_j)$ and $\pdim_{B^{\op}}(S_i^\ast) = \pdim_{B^{\op}}(S_j^\ast)$ for all $i,j$ where $S_i$ and $S_i^\ast$ denote the simple right and left $B$-modules. Then $B$ satisfies the generalized Gorenstein condition.

\begin{pro}
Since $B$ is noetherian, \cite[Theorem 8.27]{MR2455920} implies $\gldim(B) = \gldim(B^{\op})$ and thus $\pdim(S_j) = \pdim(S_1) = \pdim(S_1^\ast) = \pdim(S_j^\ast) = d$ for all $j$. In particular, it follows that $\gldim(B^{\op}) < \infty$ and from the other assumptions we get $b_1 = b_2 = \ldots = b_n = b$ for
\[
	b_j = \pdim_{B^{\op}}(S_j^\ast) = \sup\{i \mid \Ext_{B^{\op}}^i(S_j^\ast,B) \neq 0\} < \infty.
\]
Furthermore, $b = \gldim(B^{\op}) = \gldim(B) = d$. We use the fourth quadrant spectral sequence from \cite[2.2 Theorem]{levasseur_1992} as follows:
\[
	E_{p,-q}^2 = \Ext_{B^{\op}}^p\left(\Ext_B^{q}(S_j,B), B\right) \ \ \ \R \ \ \ \begin{cases} 0, & p \neq q \\ S_j, & p = q \end{cases}.
\]
In particular, $E_{p,-q}^\infty$ is nonzero precisely if $p = q$. Define
\[
	f_j = \inf\{i \mid \Ext_B^i(S_j,B) \neq 0\} < \infty.
\]
This implies all the nonzero entries of the spectral sequence can be found in the rectangle $[0,d] \times [-d,-f_j]$. 

Additionally, the $\chi$-condition guarantees that $X = \Ext_B^{f_j}(S_j,B) \neq 0$ is finite dimensional. Every finite dimensional module over $B^{\op}$ has a composition series with factors which are shifted copies of simple modules. Therefore, by short exact sequences of the first component of $\Ext_{B^{\op}}^i(-,B)$ and the definition of $b = d$ we gain $\Ext_{B^{\op}}^b(X, B) \neq 0$. Combined, this implies $E_{d,-f_j}^2 \neq 0$ which is the top right corner of the rectangle. 

The data includes the differentials $d_{p,q}^r$ of degree $(r,1-r)$ as described in \cite[1.1]{Levasseur1985}. Thus, taking homology does not change $E_{d,-f_j}^r$ for all $r$ and so $E_{d,-f_j}^2 \cong E_{d,-f_j}^\infty$ must be a nonzero entry. The spectral sequence shows now that $d = f_j$ for all $j$. In other words, we proved $d = f_1 = \ldots = f_n =: f$. 

Finally, since $d = f$ we get that $\Ext_{B}^i(S_j,B)$ is nonzero if and only if $i = d$ and $\Ext_{B}^d(S_j,B)$ is finite dimensional due to the $\chi$-condition of $B$. Repeating the argument for $B^{\op}$ shows that $\Ext_{B^{\op}}^i(S_j^\ast,B)$ is nonzero if and only if $i = d$ and $\Ext_{B^{\op}}^d(S_j^\ast,B)$ is finite dimensional. Proposition \ref{prop62} gives in this case $\Ext_B^d(S_j,B) \cong S_{m_j}^\ast$ for some $m_j \in \{1,\ldots,n\}$ for all $j$. Since
\[
	S_j = \Ext^d_{B^{\op}}(\Ext^d_B(S_j,B),B) = \Ext^d_{B^{\op}}(S_{m_j}^\ast,B),
\]
we know that $\{m_1,\ldots,m_n\} = \{1,\ldots,n\}$, and therefore $B$ satisfies the generalized Gorenstein condition. \qed
\end{pro}
\end{thm}

Since we want to apply the previous theorem to a fixed ring $B^G$ of a noetherian $\NN$-graded twisted Calabi-Yau algebra $B$ with degree zero piece $B_0 \cong k^n$, some of the assumptions are automatically satisfied. Therefore, we state the result for our purposes in the following corollary.
\begin{cor}\label{cor64}
Assume $B$ is a noetherian $\NN$-graded twisted Calabi-Yau algebra of global dimension $d \geq 1$ with degree zero piece $B_0 \cong k^n$. Suppose $\gldim(B^G) = d$ is finite and that we have the equalities $\pdim(S_i) = \pdim(S_j)$ as well as $\pdim(S_i^\ast) = \pdim(S_j^\ast)$ for all $i,j$ where $S_i = e_iB^G/(e_iB^G)_{\geq 1}$ and $S_i^\ast = B^Ge_i/(B^Ge_i)_{\geq 1}$. Then $B^G$ twisted Calabi-Yau.
\begin{pro}
Notice that $B$ and $B^{\op}$ satisfy the $\chi$-condition by Lemma \ref{lemma27}. By our convention of this section, every $g \in G$ fixes the vertices. Therefore, $(B^G)_0 = B_0 \cong k^n$. Lemma \ref{lemma214}(b) and Proposition \ref{prop217} say that $B^G$ is noetherian and that $B^G$ and $(B^G)^{\op}$ satisfy the $\chi$-condition. The result follows from Theorem \ref{thm63} applied to $B^G$ combined with Proposition \ref{prop29}.\qed
\end{pro}
\end{cor}

We finish this section by putting together all the results from Sections \ref{sectiongeneral} and \ref{section6}. By additionally requiring that all simple modules have the same projective dimension over the fixed ring we reach our goal to prove that the following three conditions of the fixed ring are equivalent.

\begin{thm}\label{thm65}
Assume $B$ is a noetherian $\NN$-graded twisted Calabi-Yau algebra of global dimension $d \geq 1$ with degree zero piece $B_0 \cong k^n$. Denote the pairwise orthogonal idempotents by $\{e_1,\ldots,e_n\}$. Let $G$ be a finite subgroup of $\Aut_{\gr}(B)$ such that $g(e_i) = e_i$ for all $g \in G$ and all $i = 1,\ldots,n$. Suppose $\pdim(S_i) = \pdim(S_j)$ and $\pdim(S_i^\ast) = \pdim(S_j^\ast)$ for the simple $B^G$-modules $S_i = e_iB^G/(e_iB^G)_{\geq 1}$ and $S_i^\ast = B^Ge_i/(B^Ge_i)_{\geq 1}$ and all $i,j = 1,\ldots,n$. Then the following are equivalent:
\begin{itemize}[topsep=0pt]
	\item[(1)]	$\gldim(B^G) = \gldim(B)$.
	\item[(2)]	$B^G$ is twisted Calabi-Yau.
	\item[(3)]	$B$ is a projective $B^G$-module.
\end{itemize}
\begin{pro}
From Proposition \ref{prop55} we know that (1) is equivalent to (3). Then the Propositions \ref{prop29} and \ref{prop52} show that (2) implies (1). The requirements of Corollary \ref{cor64} are satisfied and so (1) implies (2). \qed
\end{pro}
\end{thm}

\section{The main theorem for $\widetilde{A_{n-1}}$}\label{section7}
In this final section, we focus on the preprojective algebra $A = \Pi_Q(k)$ with $Q= \widetilde{A_{n-1}}$. We continue to use the notation $S_i = e_iA/(e_iA)_{\geq 1}$ and $S_i^\ast = Ae_i/(Ae_i)_{\geq 1}$ for the simple right and left modules of both $A$ and $A^G$ where $G$ is a finite subgroup of the graded automorphism group of $A$ such that every $g \in G$ satisfies $g(e_i) = e_i$. Another application of Theorem \ref{thm63} together with some new terminology is Proposition \ref{prop76}. This is enough to show that $A$ satisfies the required conditions for Theorem \ref{thm65}.

\begin{other}{Definition \& Remark}
Let $Q = \widetilde{A_{n-1}}$. Define $Q^G$ to be the quiver with weighted arrows which has the same vertices as $Q$ and whose arrows are the paths in $Q$ minimally generating $A^G$. Recall that $G$ is finite.

Since $G$ acts by scalar multiplication, a minimal linear combination of paths in $A = \Pi_Q(k)$ for $Q = \widetilde{A_{n-1}}$ is fixed if and only if each summand is fixed. Therefore, every generator of $A^G$ is a path and it makes sense to define $Q^G$.
\end{other}

The final part of this paper consists of proving that in case $Q = \widetilde{A_{n-1}}$ the conditions of Corollary \ref{cor64} are satisfied. In particular, we need to check that every simple module $S_i = e_iA^G/(e_iA^G)_{\geq 1}$ and $S_i^\ast = A^Ge_i/(A^Ge_i)_{\geq 1}$ has projective dimension greater than or equal to $2$. For this, we need to understand the generators and relations of $A^G$.

\begin{lemma}\label{lemma72}
Let $Q = \widetilde{A_{n-1}}$. Every vertex $e_i$ in $Q^G$ must have a purely star arrow and a purely nonstar arrow starting at $e_i$ and one of each ending at $e_i$. Further, there is a purely nonstar (purely star) path starting and ending at $e_i$ which is fixed.
\begin{pro}
We know $G = \{g_1,\ldots, g_\ell\}$ is a finite group. For $g_j$, we know $g_j(\alpha_s) = c_s^{(j)} \alpha_s$ for $s = 1,\ldots,n$ where $c_s^{(j)}$ is some root of unity. Therefore, for $j = 1,\ldots, \ell$ we have that $c_1^{(j)} c_2^{(j)} \cdots c_n^{(j)}$ is a root of unity whose order can be denoted by $n_j$. This gives that the purely nonstar path $\beta_i = (\alpha_i \alpha_{i+1} \alpha_{i+2} \cdots \alpha_n \alpha_1 \cdots \alpha_{i-1})^{n_1 \cdots n_\ell}$ starting and ending at the vertex $e_i$ must be fixed. Now, $\beta_i$ decomposes into indecomposable fixed paths in $Q^G$. The first one starts at $e_i$ and the last one ends at $e_i$. Therefore, there exists a shortest fixed purely nonstar path starting at $e_i$ as well as one ending at $e_i$ for $i = 1,\ldots,n$. Analogously, there are shortest fixed purely star paths starting at $e_i$ and ending at $e_i$.\qed
\end{pro}
\end{lemma}

\begin{cor}\label{cor73}
For $Q = \widetilde{A_{n-1}}$, the fixed ring $A^G$ has at least $3n$ generators counting the primitive orthogonal idempotents.
\begin{pro}
For every vertex $e_i$, there exist one purely nonstar and one purely star generator of $A^G$ starting at $e_i$ by Lemma \ref{lemma72} plus the idempotent $e_i$ gives three generators at $e_i$.\qed
\end{pro}
\end{cor}

Since the generators of $Q^G$ have been understood well enough we focus on the relations in $A^G$. It turns out that they have a nice form which enables us to prove the existence of a minimal relation starting and ending at each vertex. 
\begin{prop}\label{prop74}
Let $Q = \widetilde{A_{n-1}}$. Every relation in $A^G$ can be written as a sum of homogeneous (with respect to the grading in $A$) binomial relations in $A^G$.
\begin{pro}
Let $R = \beta_1 + b_2 \beta_2 + \ldots + b_m \beta_m$ be a relation in $A^G$ where the $\beta_i$ are nonzero paths in $Q^G$ and $b_2 \cdots b_m \neq 0$. We show by induction on $m$ that $R$ is a sum of binomial relations. We can understand $R$ as a relation in $A$ and therefore we can assume that all $\beta_i$ are homogeneous of the same degree in $A$, say $\ell$.

Suppose $m = 1$. This case cannot happen as a nontrivial path cannot equal zero.
If $m = 2$, then $R$ is binomial itself, $\beta_1 = -b_2 \beta_2$, so the claim follows.

Finally let $m \geq 2$. Since $R = 0$ in $A$, we can write $\beta_1 = -(b_2 \beta_2 + \ldots b_{m+1} \beta_{m+1}).$ Define $j$ and $k$ such that $\beta_1 \in e_j A e_k$. Then all $\beta_i \in e_j A e_k$ as otherwise multiplying $R$ by $e_j$ from the left and $e_k$ from the right decreases $m+1$ and induction finishes the proof. Since every path in $A$ is uniquely determined by its starting vector $e_j$, its length $\ell$ and the number of nonstar arrows, we know a basis of $(e_j A e_k)_\ell$ consists of one representative for each possible type $(0,\ell),\ (1,\ell - 1), \ldots, (\ell, 0)$. Denote the type of $\beta_i$ by $(m_i, n_i)$ for all $i$. Suppose that $(m_1,n_1) \neq (m_i,n_i)$ for all $i = 2,\ldots,m+1$. Then $\beta_1$ can be written as a linear combination of paths of different type which yields a basis of $(e_jAe_k)_{\ell}$ of size $\leq \ell$. This is a contradiction as $\dim\left((e_jAe_k)_{\ell}\right) = \ell + 1$. Hence, we can assume $(m_1,n_1) = (m_2,n_2)$. This however gives $\beta_1 = \beta_2$ and so we find the new relation $R'$ as $(b_2 - 1)\beta_2 + \ldots + b_{m+1}\beta_{m+1}$ which by induction can be written as a sum of binomial relations. We notice that since $\beta_1, \beta_2 \in A^G$, the new relation holds in $A^G$.\qed
\end{pro}
\end{prop}
This proof also shows that the minimal relations of $A^G$ are binomial of the form $\beta_1 = \beta_2$ for two paths $\beta_1, \ \beta_2 \in Q^G$. In particular, $\beta_1$ and $\beta_2$ must be of the same type in $A$. This forces $\beta_1$ and $\beta_2$ to be mixed paths.

\begin{lemma}\label{lemma75}
Let $Q = \widetilde{A_{n-1}}$. Then for every vertex $e_i$ in $A^G$ there exists a minimal relation starting at $e_i$ and a minimal relation ending at $e_i$.
\begin{pro}
Fix an idempotent $e_i \in A^G$ for some $i$. From Lemma \ref{lemma72} we obtain an indecomposable purely nonstar element $x_i \in e_iA^G$ and an indecomposable purely nonstar element $y_i \in A^Ge_i$. By construction, $x_i$ and $y_i$ are part of the set of minimal generators of $A^G$. For a vertex $e_m$, let $\nu_m \in e_mA^Ge_m$ be the shortest purely star closed path which is fixed under the action of $G$. Notice that $\nu_m$ is not necessarily indecomposable in $A^G$. Similar to the proof of Lemma \ref{lemma72}, it is easy to see that $\nu_m$ is the first power of $(\alpha_{m-1}^\ast \alpha_{m-2}^\ast \cdots \alpha_1^\ast \alpha_m^\ast \cdots \alpha_{m}^\ast)$ which is fixed. Since this construction is independent of $m$, we gain $\deg(\nu_m) = \deg(\nu_\ell)$ for all $m, \ell = 1,\ldots,n$. 

This implies that $\type(x_i \nu_{t(x_i)}) = \type(\nu_i x_i)$ and both these paths are from $e_i$ to $e_{t(x_i)}$. Thus, they are equal and the two relations arise:
\[
	x_i \nu_{t(x_i)} = \nu_i x_i \ \ \ \ \ \ \ \ \ \ \text{ and } \ \ \ \ \ \ \ \ \ \ \nu_{s(y_i)} y_i = y_i \nu_i.
\]
Since in the left relation the left hand side starts with a nonstar arrow while the right hand side starts with a star arrow, there must exist a minimal relation starting at $e_i$. Analogously, since the right relation ends with arrows of different type, there must exist a minimal relation ending at $e_i$.\qed
\end{pro}
\end{lemma}
After this setup, we are ready to prove that $A$ satisfies the missing direction from Section \ref{sectiongeneral}. The key argument established by the previous statements is that the projective dimension of each simple module is at least $2$ and therefore global dimension of $A^G$ being $2$ shows the necessary equality of all $\pdim(S_j)$ for $j = 1,\ldots,n$.
\begin{prop}\label{prop76}
Let $Q = \widetilde{A_{n-1}}$. If $\gldim(A^G) = \gldim(A) = 2$, then $A^G$ satisfies the generalized Gorenstein condition $\Ext_{A^G}^i(S_j,A^G) = \delta_{di} S_{\sigma(j)}^\ast$ for some $\sigma \in \Sym(n)$.
\begin{pro}
Denote the indecomposable projectives $e_jA^G$ by $P_j$ for $j = 1,\ldots, n$. We know the start of a minimal projective resolution of $S_j$ and with $\gldim(A^G) = 2$ it follows that
\[
\begin{tikzcd}
0 \arrow[r] & Q_2 \arrow[r] & Q_1 \arrow[r] & P_j \arrow[r] & S_j \arrow[r] & 0 
\end{tikzcd}
\]
where $Q_1$ and $Q_2$ are projective $A^G$-modules. As used before, the proof of \cite[Theorem 3.2]{MR2355031} tells us that $Q_1$ comes from the generators and $Q_2$ is determined by the relations. Lemma \ref{lemma72} shows that $A^G$ has at least two minimal generators of positive degree starting at each vertex. Hence, $Q_1 \neq 0$. By Lemma \ref{lemma75} there exists at least one minimal relation starting at each vertex. This says the kernel of the map $Q_1 \rightarrow P_j$ is nonzero or in other words, $Q_2 \neq 0$. Together, we get that $\pdim(S_j) = 2$ for all $j$. The same argument with generators and relations ending at $e_j$ works for a minimal projective resolution of $S_j^\ast$ and so $\pdim(S_j^\ast) = 2$ for all $j$. Corollary \ref{cor64} entails that $A^G$ satisfies the generalized Gorenstein condition. \qed
\end{pro}
\end{prop}

At last, we can combine the previous results to get the main theorem for $Q = \widetilde{A_{n-1}}$. We believe that the conditions are satisfied for all types $\widetilde{A_{n-1}}, \ \widetilde{D_{n-1}}$ or $\widetilde{E_{m}}$ for $m = 6,7,8$ rather than just $\widetilde{A_{n-1}}$. However, the relations of the fixed rings in the other two cases are much harder to understand. In particular, they are not binomial and it is unclear if every vertex has the needed generators and relations to get $2$ as a lower bound for the projective dimension of every simple module. Nevertheless, the following theorem motivates the decision to require $\gldim(A^G) = 2$ on the algebra side of a possible Chevalley-Shephard-Todd Theorem.

\begin{thm}\label{thm77}
For type $Q = \widetilde{A_{n-1}}$, $A = \Pi_Q(k)$ and $G$ a finite subgroup of $\Aut_{\gr}(A)$ such that $g(e_i) = e_i$ for all $g \in G$ and all $i = 1,\ldots,n$, the following conditions are equivalent:
\begin{itemize}[topsep=0pt]
	\item[(1)]	$\gldim(A^G) = \gldim(A)$ (i.e. $\gldim(A^G) = 2$).
	\item[(2)]	$A^G$ is twisted Calabi-Yau.
	\item[(3)]	$A$ is a projective $A^G$-module.
\end{itemize}
\begin{pro}
Proposition \ref{prop76} shows that (1) implies (2) using Proposition \ref{prop29}. The other implications hold due to \ref{prop52} and \ref{prop55}.\qed
\end{pro}
\end{thm}

\bibliographystyle{alpha}
\bibliography{POTFROAPAbib}

\end{document}